\documentclass[a4paper,twoside]{article}
\usepackage{amssymb}
\usepackage{amsmath}
\usepackage{amsfonts}
\usepackage{}
\usepackage{color}
\usepackage{mathrsfs}
\usepackage{amsmath,amsfonts,amssymb,amsthm}
\usepackage{enumitem}
\usepackage[bookmarks]{hyperref}

\usepackage[center]{titlesec}
\usepackage{tocbibind}
\usepackage{geometry}
\numberwithin{equation}{section}

\newtheorem{thm}{Theorem}[section]
\newtheorem{cj}{Conjecture}[section]
\newtheorem{defi}[thm]{Definition}
\newtheorem{lem}[thm]{Lemma}
\newtheorem{cor}[thm]{Corollary}
\newtheorem{rem}[thm]{Remark}
\newtheorem{prop}[thm]{Proposition}

\def\bb{\mathbb}
\def\ca{\mathcal}

\def\e{\epsilon}

\def\scr{\mathscr}
\def\bbd{\boldsymbol}

\def\dirac{\boldsymbol{\delta}}

\newcommand{\bra}[1]{\left(#1\right)}

\newcommand{\abs}[1]{\left\vert#1\right\vert}
\newcommand{\set}[1]{\left\{#1\right\}}

\def\bthm{\begin{thm}\def\ethm{\end{thm}}}
\def\brem{\begin{rem}\def\erem{\end{rem}}}

\def\bdefi{\begin{defi}\sl{}\def\edefi{\end{defi}}}
\def\blem{\begin{lem}\def\elem{\end{lem}}}

\def\bprop{\begin{prop}\def\eprop{\end{prop}}}

\def\bitm{\begin{itemize}} \def\eitm{\end{itemize}}
\def\benu{\begin{enumerate}} \def\eenu{\end{enumerate}}

\def\bpf{\begin{proof}}\def\epf{\end{proof}}
\def\beq{\begin{equation}}\def\eeq{\end{equation}}
\def\beqs{\begin{eqnarray}}\def\eeqs{\end{eqnarray}}
\def\beqsnl{\begin{eqnarray*}}\def\eeqsnl{\end{eqnarray*}}

\usepackage{fancyhdr}
\begin{document}
\title{
The Multifractal Nature of Boltzmann Processes}

\footnotetext{This work was supported by China Scholarship Council.}
\author{Liping Xu}
%\date{}%(空着表示不显示日期)
\maketitle

\begin{abstract}
We consider the spatially homogeneous Boltzmann equation for (true) hard
and moderately soft potentials.
We study the pathwise properties of the stochastic process $(V_t)_{t\geq 0}$, which describes the time
evolution of the velocity of a typical particle. We show that this process is almost surely multifractal
and we compute its spectrum of singularities.
For hard potentials, we also compute the multifractal spectrum of the position process $(X_t)_{t\ge0}$.
\end{abstract}

\section{Introduction}
\quad~~The Boltzmann equation is the main model of kinetic theory. It describes the time evolution
of the density $f_t(x,v)$ of particles with position $x\in\bb{R}^3$ and velocity $v\in \bb{R}^3$ at time $t\ge 0$,
in a gas of particles interacting through binary collisions.
In the special case where the gas is initially spatially homogeneous, this property propagates with time,
and $f_t(x,v)$ does not depend on $x$. We refer to the books by Cercignani \cite{cercignani1988boltzmann} and Villani \cite{villani2002review} for many details on the physical and mathematical theory of this equation, see also the review paper by Alexandre \cite{alexandre2009review}.

Tanaka gave in \cite{tanaka1978probabilistic} a probabilistic
interpretation of the case of Maxwellian molecules:
he constructed a Markov process $(V_t)_{t\ge 0}$, solution to a Poisson-driven stochastic differential equation, and such that the law of $V_t$ is $f_t$ for all $t\geq 0$. Such a process
$(V_t)_{t\ge 0}$ has a richer structure than the Boltzmann equation, since it contains some information on the history of particles. Physically, $(V_t)_{t\ge 0}$ is interpreted as the time-evolution of the velocity of a typical particle. Fournier and M\'{e}l\'{e}ard \cite{fournier2001markov}
extended Tanaka's work to non-Maxwellian molecules, see the last part of paper by Fournier \cite{fournier2012finiteness} for up-to-date results.

In the case of long-range interactions, that is when particles interact through a repulsive force in $1/r^s$ (for some $s>2$), the Boltzmann equation presents a singular integral (case without cutoff). The reason is that the corresponding process
$(V_t)_{t\ge 0}$ jumps infinitely often, i.e. the particle is subjected to infinitely many collisions,
on each time interval. In some sense, it behaves, roughly, like a L\'{e}vy process.

The H\"{o}lder regularity of the sample paths of stochastic processes was first studied by Orey and Taylor \cite{orey1974often} and
Perkins \cite{perkins1983hausdorff}, who showed that the fast and slow points of Brownian motion are located on random sets of times, and they showed that the sets of points with a given pointwise regularity have a fractal nature. Jaffard \cite{jaffard1999multifractal}
showed that the sample paths of most L\'{e}vy processes are multifractal functions and he obtained their spectrum of singularities. This spectrum is almost surely  deterministic: of course, the sets with a given pointwise regularity are extremely complicated, but their Hausdorff dimension is deterministic. Let us also mention the article by Balan\c{c}a
\cite{balancca2014fine}, in which he extended the results
(and simplified some proofs) of Jaffard \cite{jaffard1999multifractal}.

What we expect here is that $(V_t)_{t\ge 0}$ should have the same spectrum as a well-chosen L\'{e}vy process.
This is of course very natural (having a look at the shape of the jumping SDE satisfied by $(V_t)_{t\ge 0}$).
There are however many complications, compared to the case of L\'{e}vy processes, since we loose all the independence and stationarity
properties that simplify many computations and arguments.
We will also compute the multifractal spectrum of the position process $(X_t)_{t\ge0}$, defined by
$X_t=\int_0^tV_sds$, which appears to have multifractal sample paths as well.

By the way, let us mention that, though there are many papers computing the multifractal spectrum of some quite complicated objects, we are not aware of any work concerning general Markov processes, that is, roughly, solutions to jumping (or even non jumping) SDEs. In this paper, we study the important case of the Boltzmann process, as a physical example of jumping SDE. Of course, a number of difficulties have to be encompassed, since the model is rather complicated.
However, we follow, adapting everywhere to our situation, the main ideas of Jaffard \cite{jaffard1999multifractal} and Balan\c{c}a \cite{balancca2014fine}.

Let us finally mention that Barral, Fournier, Jaffard and Seuret \cite{barral2010pure} studied a very specific ad-hoc Markov process,  showing  that quite simple processes may have a random spectrum that depends heavily on the values taken by the process.

\subsection{The Boltzmann equation}
\quad~~We consider a 3-dimensional spatially homogeneous  Boltzmann equation, which depicts the density $f_t(v)$ of particles in a gas, moving with velocity $v\in\bb{R}^3$ at time $t\geq 0$. The density $f_t(v)$ solves  \beq\label{Bol}
\partial_tf_t(v)=\int_{\bb{R}^3}dv_{\ast}\int_{\bb{S}^2}d\sigma B(|v-v_{\ast}|,\cos\theta)[f_{t}(v^{\prime})f_{t}(v_{\ast}^{\prime})-f_t(v)f_t(v_{\ast})],
\eeq
where
\beq\label{Bol1}
v^{\prime}=\frac{v+v_{\ast}}{2}+\frac{|v-v_{*}|}{2}\sigma,~v_{\ast}^{\prime}
=\frac{v+v_{\ast}}{2}-\frac{|v-v_{\ast}|}{2}\sigma,~\textrm{and}~ \cos\theta=\Big\langle\frac{v-v_{\ast}}{|v-v_{\ast}|},\sigma \Big\rangle.
\eeq

The cross section $B(|v-v_{\ast}|,\cos\theta)\geq 0$ depends on the type of interaction between particles. It only depends on $|v-v_{\ast}|$ and on the cosine of the deviation angle $\theta$.~
Conservations of mass,~momentum and kinetic energy hold for reasonable solutions and we may assume without loss of generality that $\int_{\mathbb{R}^3}f_{0}(v)dv=1$.
We will assume that there is a measurable function $\beta:(0,\pi]\rightarrow\bb{R}_+$ such that
\beq\label{con}
\left\{\begin{array}{l} B(|v-v_{\ast}|,\cos\theta)\sin\theta={|v-v_{\ast}|}^{\gamma}\beta(\theta),\\
\exists~0<c_0<C_0,~\forall~\theta\in(0,\pi/2], ~c_0\theta^{-1-\nu}\leq\beta(\theta)\leq C_0\theta^{-1-\nu},\\
\forall~\theta\in(\pi/2,\pi),~\beta(\theta)=0,
\end{array}
\right.\eeq
 for some $\nu\in(0,1)$, and $\gamma\in(-1,1)$ satisfying $\gamma+\nu>0$. The last assumption on the function $\beta$ is not a restriction and can be obtained by  symmetry as in \cite{alexandre2009review}.
Note that, when particles collide by pairs due to a repulsive force proportional to $1/r^s$ for some $s>2$, assumption (\ref{con}) holds with $\gamma=(s-5)/(s-1)$ and $\nu=2/(s-1)$.
Here we will be focused on the cases of hard potentials $(s>5)$, Maxwell molecules $(s=5)$ and moderately soft potentials $(3<s<5)$.

Next, we give the definition of weak solutions of (\ref{Bol}). We define $\mathcal {P}_p(\bb{R}^3)$ as the set of all probability measures $f$ on $\bb{R}^3$ such that $m_p(f):=\int_{\bb{R}^3}|v|^pf(dv)<\infty$.
\bdefi
Assume \emph{(\ref{con})} is true for some $\nu\in(0,1), \gamma\in(-1,1)$.
 A measurable family of probability measures $(f_t)_{t\geq 0}$ on $\bb{R}^3$ is called a \textit{weak solution} to \emph{(\ref{Bol})} if it satisfies the following two conditions.
\begin{itemize}
  \item For all $t\geq 0$,
  \beq
\int_{\bb{R}^3}vf_t(dv)=\int_{\bb{R}^3}vf_0(dv)~~\text{and}~~
\int_{\bb{R}^3}|v|^2f_t(dv)=\int_{\bb{R}^3}|v|^2f_0(dv)<\infty.
\eeq
  \item For any bounded globally Lipschitz-continuous function $\phi:\bb{R}^3\rightarrow\bb{R}$, any $t\geq 0$,
  \beq
  \int_{\bb{R}^3}\phi(v)f_t(dv)=\int_{\bb{R}^3}\phi(v)f_0(dv)+
  \int_0^t\int_{\bb{R}^3}\int_{\bb{R}^3}L_B\phi(v,v_\ast)f_s(dv_\ast)f_s(dv)ds,
  \eeq
  where $v^\prime$ and $\theta$ are defined by \emph{(\ref{Bol1})}, and
  \[L_B\phi(v,v_\ast):=\int_{\mathbb{S}^2}B(|v-v_\ast|,\cos\theta)
  (\phi(v^\prime)-\phi(v))d\sigma.\]
\end{itemize}
\edefi

The existence of a weak solution to (\ref{Bol}) is now well established (see \cite{villani1998new}
and \cite{lu2012measure}).
In particular, when $\gamma \in (0,1)$, it is shown in \cite{lu2012measure} that for any $f_0 \in
\mathcal {P}_2({\bb{R}^3})$, there exists a weak solution $(f_t)_{t\geq 0}$ to (\ref{Bol}) satisfying
$\sup_{t\ge t_0}m_p(f_t)<\infty$ for all $p\ge2$, all $t_0>0$.
Some uniqueness results can be found in \cite{fm09}.

\subsection{The Boltzmann process}
\quad~~We first parameterize (\ref{Bol1}) as in \cite{fournier2002stochastic}. For each $x \in\bb{R}^{3}\setminus{\{0\}}$,
we consider the vector $I(x)\in\bb{R}^3$ such that $|I(x)|=|x|$ and $I(x)\perp x$. We also set $J(x)=\frac{x}{|x|}\wedge I(x)$. The triplet $(\frac{x}{|x|},\frac{I(x)}{|x|},\frac{J(x)}{|x|})$ is an orthonormal basis of $\bb{R}^3$. Then for $x,v,v_\ast\in\bb{R}^3$,
$\theta\in[0,\pi),~\varphi\in[0,2\pi)$, we set
\beq\label{para}
\left\{\begin{array}{l} \Gamma(x,\varphi):=(\cos\varphi) I(x)+(\sin\varphi) J(x),\\
v^\prime(v,v_\ast,\theta,\varphi):=v-\frac{1-\cos\theta}{2}(v-v_\ast)
+\frac{\sin\theta}{2}\Gamma(v-v_\ast,\varphi),\\
{\bf a}(v,v_\ast,\theta,\varphi):=v^\prime(v,v_\ast,\theta,\varphi)-v.
\end{array}
\right.\eeq
Let us observe at once that $\Gamma(x,\varphi)$ is orthogonal to $x$ and has the same norm as $x$,
from which it is easy to check that
\beq\label{num}
|{\bf a}(v,v_*,\theta,\varphi)|=\sqrt{\frac{1-\cos\theta}{2}}|v-v_*|.
\eeq
\bdefi\label{def}
Let $(f_t)_{t\ge0}$ be a weak solution to \emph{(\ref{Bol})}. On some probability space
$(\Omega,\mathcal {F},(\mathcal {F}_t)_{t\geq0}, \mathbb{P})$, we consider
a $\mathcal {F}_0$-measurable random variable $V_0$ with law $f_0$, a Poisson measure $N(ds,dv,d\theta,d\varphi,du)$ on $[0,\infty)\times\bb{R}^3\times(0,\pi/2]\times[0,2\pi)\times[0,\infty)$ with intensity $dsf_s(dv)\beta(\theta)d\theta d\varphi du$.
A c\`adl\`ag $(\mathcal {F}_t)_{t\geq0}$-adapted process $(V_t)_{t\geq 0}$ with values in $\bb{R}^3$ is then called a \textit{Boltzmann process} if it solves
\beq\label{Bol3}
V_t=V_0+\int_{0}^{t}\int_{\bb{R}^3}\int_0^{\pi/2}\int_0^{2\pi}
\int_0^{\infty}{\bf a}(V_{s-},v,\theta,\varphi)\bbd{1}_{\{u\leq|V_{s-}-v|^\gamma\}}
N(ds,dv,d\theta,d\varphi,du).
\eeq
\edefi

From Proposition 5.1 in \cite{fournier2012finiteness}, we have slightly different results for different potentials: when $\gamma\in(0,1)$, i.e. hard potentials,
we can associate a Boltzmann process to any weak solution to {(\ref{Bol})}, but when $\gamma\in(-1,0)$,
i.e. moderately soft potentials, we can only prove existence of  a weak solution  to {(\ref{Bol})} to which it is possible to associate a Boltzmann process.

 \bprop\label{prop}
 Let $f_0\in\mathcal {P}_2({\bb{R}^3})$. Assume \emph{(\ref{con})} for some $\gamma\in(-1,1),~\nu\in(0,1)$.
 \begin{itemize}
   \item If $\gamma\in(0,1)$, for any weak solution $(f_t)_t\geq0$ to \emph{(\ref{Bol})} starting from $f_0$
   and satisfying
   $$\mbox{ for all $p\ge2$, all $t_0>0$, } \ \ \sup_{t\ge t_0}m_p(f_t)<\infty,$$
   there exist  a  probability space $(\Omega,\mathcal {F},(\mathcal {F}_t)_{t\geq0},\mathbb{P})$,
   a $(\mathcal {F}_t)_{t>0}$-Poisson measure $N(ds,dv,d\theta,d\varphi,du)$ on $[0,\infty)\times\bb{R}^3\times(0,\pi/2]\times[0,2\pi)\times[0,\infty)$
   with intensity $dsf_s(dv)\beta(\theta)d\theta d\varphi du$ and a c\`adl\`ag $(\mathcal {F}_t)_{t\geq0}$-adapted
   process $(V_t)_{t\geq 0}$ satisfying $\mathcal {L}(V_t)=f_t$ for all $t\geq0$ and solving \emph{(\ref{Bol3})}.

    \item If $\gamma\in(-1,0]$, assume additionally that $f_0\in\mathcal {P}_p(\bb{R}^3)$ for some $p>2$.
    There exist a probability space, a Poisson measure $N$ and a c\`adl\`ag adapted process $(V_t)_{t\geq 0}$
    as in the previous case, satisfying $\mathcal {L}(V_t)=f_t$ for all
    $t\geq0$ and solving \emph{(\ref{Bol3})}.
 \end{itemize}
 \eprop

The Boltzmann equation depicts the velocity distribution of a dilute gas which is made up of a large number of molecules.
So, the corresponding Boltzmann process $(V_t)_{t\geq0}$ represents the time evolution of the velocity of a typical particle.
When this particle collides with another one, its velocity changes suddenly. It is thus a jump process.
\subsection{Recalls on multifractal analysis}
\quad~~In this part, we   recall the definition of the main  objects in multifractal analysis.
\bdefi\label{ddhe}
A locally bounded function $g:[0,1]\rightarrow\bb{R}^3$ is said to belong to the
\textit{pointwise H\"{o}lder space} $C^{\alpha}(t_0)$ with $t_0\in[0,1]$
and $\alpha\notin\bb{N}$, if there exist $C>0$ and a polynomial $P_{t_0}$ of degree less than $\lfloor \alpha\rfloor$, such that for some neighborhood $I_{t_0}$ of $t_0$,
\[|g(t)-P_{t_0}(t)|\leq C|t-t_0|^\alpha, ~\forall~t\in I_{t_0}.\]
The \textit{pointwise H\"{o}lder exponent} of $g$ at point $t_0$ is given by
\[h_g(t_0)=\sup\{\alpha>0:g\in C^\alpha(t_0)\},\]
 where by convention $\sup\emptyset=0$. The level sets of the pointwise H\"{o}lder exponent of the function $g$ are called the \textit{iso-H\"{o}lder sets} of $g$, and are denoted,
 for any $h\ge0$, by  \[E_g(h)=\{t\geq0:h_g(t)=h\}.\]
\edefi

We now recall the definition of the Hausdorff measures and dimension, see \cite{falconer2007fractal}
for   details.
\bdefi\label{def1}
Given a subset $A$ of $\bb{R}$, given $s>0$ and $\epsilon>0$, the $s$-Hausdorff pre-measure $\ca{H}_{\epsilon}^{s}$ using balls of radius less than $\epsilon$ is given by
\[
\ca{H}_{\epsilon}^{s}(A)=\inf\left\{\sum_{i\in J}|I_i|^s:~(I_i)_{i\in J}\in\mathscr{P}_\epsilon(A)\right\},
\]
where $\mathscr{P}_\epsilon(A)$ is the set of all countable coverings of $A$ by intervals with length at most $\epsilon$.
The $s$-Hausdorff measure of $A$ is defined by \[\ca{H}^s(A)=\lim_{\epsilon\rightarrow0}\ca{H}_{\epsilon}^{s}(A).\]
Finally the \textit{Hausdorff dimension} of $A$ is defined by \[\dim_H(A):=\inf\{s\geq0:\ca{H}^s(A)=0\}=\sup\{s\geq0:\ca{H}^s(A)=+\infty\},\] and by convention $\dim_H\emptyset=-\infty$.
\edefi
We   use the concept of spectrum of singularities to describe the distribution of the singularities of a function $g$.
\bdefi
Let $g:[0,1]\rightarrow\bb{R}^3$ be a locally bounded function.
The \textit{spectrum of singularities} (or {\em multifractal spectrum)} of $g$ is the function $D_g:\bb{R}_+\rightarrow\bb{R}_+\cup{\{-\infty\}}$ defined by   \[D_g(h)=\dim_H(E_g(h)).\]
\edefi

The iso-H\"{o}lder sets $E_g(h)$ are random for most studied stochastic processes,
but almost always have an a.s. deterministic Hausdorff dimension,
as in the case of L\'{e}vy processes  \cite{jaffard1999multifractal}.

\subsection{Main Results}
\quad~~Now, we give the main results of this paper.
\bthm\label{m}
We assume \emph{(\ref{con})} for some $\gamma\in(-1,1)$, some $\nu\in(0,1)$ with $\gamma+\nu>0$.
We consider some initial condition $f_0\in\mathcal {P}_2(\bb{R}^3)$ and assume that it is not a Dirac mass. If $\gamma\in(-1,0]$,
we moreover assume that $f_0\in\mathcal {P}_p(\bb{R}^3)$ for some $p>2$.
We consider a Boltzmann process $(V_t)_{t\in [0,1]}$ as introduced in \emph{Proposition \ref{prop}}. Almost surely,
for all $h\geq 0$,
\beq
D_V(h)=        \left\{\begin{array}{ll}
           \nu h & \text{if}\quad 0\leq h\leq 1/\nu,  \\
           -\infty & \text{if} \quad h> 1/\nu.
           \end{array}
         \right.
\eeq
\ethm

The condition that $f_0$ is not a Dirac mass is important: if $V_0=v_0$ a.s. for some
deterministic $v_0 \in \bb{R}^3$, then $V_t=v_0$ for all $t\geq 0$ a.s. (which is a.s. a $C^\infty$ function on
$[0,\infty)$).

It is obvious from the proof that the spectrum of singularities is {\it homogeneous}:
we could prove similarly that a.s., for any $0\leq t_0<t_1<\infty$, all $h\geq 0$,
$\dim_H(E_V(h)\cap[t_0,t_1])=D_V(h)$.

Finally, it is likely that the same result holds true for very soft potentials. However,
there are several technical difficulties, and the proof would be much more intricate.

Now we exhibit the multifractal spectrum of the position process. For simplicity, we only consider the case of hard potentials.

\bthm\label{pos} We assume \emph{(\ref{con})} for some $\gamma\in(0,1)$ and some $\nu\in(0,1)$.
We consider some initial condition $f_0\in\mathcal {P}_2(\bb{R}^3)$ and assume that it is not a Dirac mass.
We consider a Boltzmann process $(V_t)_{t\in[0,1]}$ as introduced in \emph{Proposition \ref{prop}} and
introduce the associated position process $(X_t)_{t\in [0,1]}$ defined by
$\displaystyle X_t=\int_0^tV_s ds$. Almost surely, for all $h\ge0$,
\beq
D_X(h)=        \left\{\begin{array}{ll}
           \nu(h-1) & \text{if}\quad 1\leq h\leq \frac{1}{\nu}+1,  \\
           -\infty & \text{if} \quad h> \frac{1}{\nu}+1\quad\text{or}\quad 0\le h<1.
           \end{array}
         \right.
\eeq
\ethm

This result is very natural once Theorem \ref{m} is checked: we expect that at some given time $t$,
the pointwise exponent of $X$ is the one of $V$ plus 1.
However, this is not always true: for instance, as
can be seen on the simple example of the {\it chirp} function $g(x)=x\sin(1/x)$: its pointwise exponent at 0 is 1,
while its primitive has a pointwise exponent equal to 3 at 0.
Balan\c{c}a \cite{balancca2014fine} has shown that such an oscillatory phenomenon may occur for L\'evy processes,
but on a very small set of points.

\bdefi
Let $g:[0,1]\to\mathbb{R}^3$ be a locally bounded function and let $G(t)=\int_0^t g(s)ds$.
For all $h\ge0$, we introduce the sets
\begin{equation}
\label{defEosc}
E^{cusp}_{g}(h)=\{t\in E_{g}(h):h_{G}(t)=1+h_g(t)\} \ \mbox{ and }  \ E^{osc}_{g}(h)=\{t\in E_{g}(h):h_{G}(t)>1+h_g(t)\}.
\end{equation}
\edefi

The times $t\in E^{cusp}_{g}(h) $ are referred to as cusp singularities, while the times $t\in E^{osc}_{g}(h) $
are called oscillating singularities.
Observe that $E_{g}(h)=E^{cusp}_{g}(h)\cup E^{osc}_{g}(h)$, the union being disjoint: this follows from the fact that
obviously, for all $t\in [0,1]$, $h_G(t) \geq h_g(t)+1$. We will prove the following.

\bthm\label{pos2} Under the assumptions of Theorem \ref{pos}, we have almost surely:
\begin{itemize}
\item
  for all $h\in [1/(2\nu),1/\nu)$, $\dim_H\Big(E^{osc}_{V}(h)\Big) \le 2h\nu-1$,
\item
 for all $h\in[0,1/(2\nu)) \cup (1/\nu,+\infty]$, $E^{osc}_{V}(h)=\emptyset$,
\item
  for all $h\in [0,1/\nu]$, $\dim_H\Big(E^{cusp}_{V}(h)\Big) =  h\nu $.
 \end{itemize}
 \ethm

Actually, we will first prove Theorem \ref{pos2} which, together with Theorem \ref{m}, implies Theorem \ref{pos}.

\section{Localization of the problem}

\quad~~In the following sections,
we consider a Boltzmann process $(V_t)_{t\in [0,1]}$ associated to a weak solution $(f_t)_{t\in[0,1]}$ to (\ref{Bol}),
and driven by a Poisson measure $N(ds,dv,d\theta,d\varphi,du)$ on $[0,1]\times\bb{R}^3\times(0,\pi/2]\times[0,2\pi)\times[0,\infty)$
with intensity $dsf_s(dv)\beta(\theta)d\theta d\varphi du$.

For $B\ge1$, setting $H_B(v)=\frac{|v|\wedge B}{|v|}v$, we define, for $t\in [0,1]$,
\beq\label{newprocess}
V_t^B:=V_0+\int_0^t\int_{\bb{R}^3}\int_0^{\pi/2}\int_0^{2\pi}\int_0^{\infty}
{\bf a}(H_B(V_{s-}),v,\theta,\varphi)\bbd{1}_{\{u\le|H_B(V_{s-})-v|^\gamma\}}N(ds, dv,d\theta,d\varphi,du),
\eeq
where ${\bf a}$ is defined in (\ref{para}). We define the corresponding position process, for $t\in [0,1]$, as
\beq\label{pos-proce}
X_t^B=\int_0^tV_s^Bds.
\eeq

In the rest of the paper, we will check the following two localized claims.

\bprop\label{loc}
Let $B\ge1$ be fixed. We assume \emph{(\ref{con})} for some $\gamma\in(-1,1)$, some $\nu\in(0,1)$ with $\gamma+\nu>0$. We consider the localized process introduced in (\ref{newprocess}). Almost surely, for all $h\ge0$,
\beq
D_{V^B}(h)=        \left\{\begin{array}{ll}
           \nu h & \text{if}\quad 0\leq h\leq 1/\nu,
           \nonumber \\
           -\infty & \text{if} \quad h>1/\nu.
           \end{array}
         \right.
\eeq
\eprop
\bprop\label{locpo}
Let $B\ge1$ be fixed. We assume \emph{(\ref{con})} for some $\gamma\in(0,1)$, some $\nu\in(0,1)$. We consider the
localized process $(V_t^B)_{t\ge0}$ defined in (\ref{newprocess}). Then almost surely,
\begin{itemize}
\item
  for all $h\in [1/(2\nu),1/\nu)$, $\dim_H\Big(E^{osc}_{V^B}(h)\Big) \le 2h\nu-1$,
\item
 for all $h\in[0,1/(2\nu)) \cup (1/\nu,+\infty]$, $E^{osc}_{V^B}(h)=\emptyset$,
\item
  for all $h\in [0,1/\nu]$, $\dim_H\Big(E^{cusp}_{V^B}(h)\Big) =  h\nu $.
 \end{itemize}
\eprop

Once these propositions are verified, Theorems \ref{m} and \ref{pos2} are
immediately deduced.
\bpf[Proof of Theorems \ref{m} and \ref{pos2}]
Since $\sup_{[0,1]} |V_t| <+\infty$ a.s. (because $V$ is a c\`adl\`ag process), the event
$\Omega_B=\{\sup_{[0,1]} |V_t|  \leq B\}$ a.s. increases to $\Omega$ as $B$ increases to infinity.
But on $\Omega_B$, we obviously have that $(V_t^B)_{t \in [0,1]}=(V_t)_{t\in [0,1]}$.
Hence on $\Omega_B$, it holds that for all $h \in [0,+\infty]$,
$D_{V}(h)=D_{V^B}(h)$, $\dim_H(E^{osc}_{V}(h))=\dim_H(E^{osc}_{V^B}(h))$
and $\dim_H(E^{cusp}_{V}(h))=\dim_H(E^{cusp}_{V^B}(h))$.
The conclusion then follows from the above two propositions.
\epf

We thus fix $B\ge1$ for the rest of the paper.

\section{Study of the velocity process}

\subsection{Preliminary}
\quad~~First, we need to bound $f_t$ from below.
\blem\label{le}
There exist $a,b,c >0$, such that for any $w\in \bb{R}^3$, any $t\in[0,1]$,
\beq
f_t(\scr{H}_{w})\geq b,
\eeq
where $\scr{H}_{w}=\{v\in \bb{R}^3:|v-w|\geq a, |v|\le c\}$.
\elem
\bpf
As $f_0$ is not a Dirac mass, there exist $v_1\neq v_2$ such that $v_1, v_2\in \textrm{Supp}f_0$. We set $a=\frac{|v_1-v_2|}{6}$.

\vskip0.2cm\emph{Step 1.} We first show that there exists $b>0$, such that for all $w\in\bb{R}^3$, $t\in [0,1]$, $f_t(\{v:|v-w|\geq a\})\geq 2b$.
  First, if $|w|\ge \sqrt{2m_2(f_0)}+a=:M$, recalling that $m_2(f_t)=m_2(f_0)$ for all $t\ge0$,
\begin{align*}
f_t(\{v:|v-w|\geq a\})\ge&\ f_t(\{v:|v|\le |w|-a\})= \ 1-f_t(\{v:|v|>|w|-a\})\\
\ge&\ 1-\frac{m_2(f_0)}{(|w|-a)^2} \geq 1-\frac{m_2(f_0)}{2m_2(f_0)} = \frac{1}{2}.
\end{align*}
Next, we consider a bounded nonnegative globally Lipschitz-continuous function $\phi:\bb{R}_+\rightarrow[0,1]$,
such that for all $v>0$, $\bbd{1}_{B(0,a)^c}(v)\geq \phi(|v|)\geq\bbd{1}_{B(0,2a)^c}(v)$,
and define $F(t,w)=\int_{\bb{R}^3}\phi(|w-v|)f_t(dv)$.
We know that $t\mapsto F(t,w)$ is continuous for each $w\in\bb{R}^3$ by Lemma 3.3 in \cite{fournier2012finiteness}.
Moreover, $F(t,w)$ is (uniformly in $t$) continuous in $w$ by the Lipschitz-continuity of $\phi$. So $F(t,w)$ is continuous on $[0,1]\times\bb{R}^3$.
Since for all $t>0$, $\textrm{Supp}f_t=\bb{R}^3$ by Theorem 1.2 in \cite{fournier2012finiteness}, we get $F(t,w)\geq f_t(B(w,2a)^c)>0$,
$\forall~(t,w)\in(0,1]\times \overline{B(0,M)}$.
When $t=0$, recalling that $v_1, v_2\in \textrm{Supp}f_0$ and $a=\frac{|v_1-v_2|}{6}$, we easily see that for all $w\in\bb{R}^3$,
either $B(v_1,a)\subset B(w,2a)^c$ or $B(v_2,a)\subset B(w,2a)^c$,
whence $F(0,w)\ge \min\{f_0(B(v_1,a)),f_0(B(v_2,a))\}>0$. Since $[0,1]\times \overline{B(0,M)}$ is compact and $F(t,w)$ is continuous,
there exists $b_1>0$, such that $f_t(B(w,a)^c)\geq F(t,w)\geq b_1$ for all $(t,w)\in [0,1]\times \overline{B(0,M)}$.
So we conclude by choosing $b=\min(\frac{1}{2}, b_1)/2$.

\vskip0.2cm\emph{Step 2.} We now conclude. Using Step 1,
\begin{align*}
f_t(\{v:|v-w|\geq a,|v|\leq c\})
\ge&\ f_t(\{v:|v-w|\geq a\})-f_t(\{v:|v|>c\})
\ge\ 2b-\frac{m_2(f_0)}{c^2}.
\end{align*}
So, we complete the proof by taking $c= \sqrt{\frac{m_2(f_0)}{b}}$.
\epf

\subsection{random fractal sets associated with the Poisson process}
\quad~~
First, we introduce some notations. Recall that $h_{V^B}, E_{V^B}, D_{V^B}$ respectively the H\"{o}lder exponent,
iso-H\"{o}lder set and spectrum of singularities of the Boltzmann process $(V_t^B)_{t\in [0,1]}$.
The notation $\scr{L}$ represents the Lebesgue measure.
$\mathcal {J}$ designates the set of the jump times of the process $V^B$, that is,
\[\mathcal {J}:=\{s\in[0,1]:|\Delta V_s^B|\neq0\}.\] For $m\ge1$, we also introduce
\[\ca{J}_m:=\{s\in\ca{J}:|\Delta V_s^B|\le2^{-m}\},~\widetilde{\ca{J}}_m:=\{s\in\ca{J}:2^{-m-1}<|\Delta V_s^B|\le2^{-m}\}.\]
For $\delta>0$ and $m\ge1$, we define the sets
 \[A_\delta^m:=\bigcup_{s\in \ca {J}_m}[s-|\Delta V_s^B|^{\delta},s+|\Delta V_s^B|^{\delta}],~
\widetilde{A}_\delta^m:=\bigcup_{s\in\widetilde{\ca{J}}_m}[s-|\Delta V_s^B|^\delta, s+|\Delta V_s^B|^\delta].\]
Finally, for $\delta>0$, we define
\beq\label{def2}A_{\delta}=\limsup_{m\rightarrow+\infty}A_\delta^m
=\limsup_{m\rightarrow+\infty}\widetilde{A}_\delta^m.\eeq
The main result of this subsection states that
\bprop\label{rfs}
We have a.s. the following properties:
\begin{enumerate}[label=\emph{(\arabic*)}]
\item for all $\delta\in(0,\nu), A_\delta\supset[0,1]$,
\item there exists a (random) positive sequence $(\epsilon_m)_{m\ge1}$ decreasing to 0, such that
\[\scr{L}\left(A_\nu^*~\bigcap~[0,1]\right)=1,\] where we use the notation $A_\delta^*=\limsup_{m\rightarrow +\infty} \widetilde{A}_{\delta(1-\epsilon_m)}^m$, for all $\delta\in(0,\infty)$.
\end{enumerate}
\eprop
\brem\label{rema}
We observe at once that for any $\delta>\delta^\prime>0$,
$A_\delta\subset A_\delta^*\subset A_{\delta^\prime}$.
\erem
The reason why we study $A_\delta$ comes from  the following heuristics:
if $t\in A_\delta$ with $\delta$ large, then $t$ is rather close to
many large jump times of $V^B$, so that $V^B$ will not be very regular at $t$.
On the contrary, if $t$ does only belong to those $A_\delta$'s with $\delta$ small, this means that
$t$ is rather far away from the jumps of $V^B$, so that $V^B$ will be rather regular at $t$.

We introduce $A_\delta^*$ (which resembles very much $A_\delta$) for technical reasons, mainly because
at the {\it critical} value $\delta=\nu$, we cannot prove (and it may be false)
that $A_\nu$ has a full Lebesgue measure.

The rest of this subsection is devoted to proving this proposition.
We first recall the Shepp lemma, first discovered in \cite{shepp1972covering},
in the version used in \cite{jaffard1999multifractal}.
\blem\label{shep}
 We consider a Poisson measure $\pi(ds,dy)=\sum_{s\in \mathscr{D}}\dirac_{(s,y_s)}$ on $[0,1]\times(0,1)$ with
intensity $ds\mu(dy)$, where $\mu$ is a measure on $(0,1)$.
 We consider the set $U=\cup_{s\in \mathscr{D}}(s-y_s,s+y_s)$.
 If \[\int_0^1\exp\left({2\int_t^1\mu((y,1))dy}\right)dt=+\infty,\]
then almost surely, $[0,1]\subset U$.
\elem
We write $N=\sum_{s\in \mathscr{D}}\dirac_{(s,v_s,\theta_s,\varphi_s,u_s)}$,
where $v_s,\theta_s,\varphi_s,u_s$ are the quanta corresponding to the jump time $s\in\mathscr{D}$.
For convenience, we consider this Poisson measure by adding a family of independent variables $(x_s)_{s\in\mathscr{D}}$,
which are uniformly distributed in $[0,1]$ and independent of $v_s,\theta_s,\varphi_s,u_s$,
so that $O:=\sum_{s\in \mathscr{D}}\dirac_{(s,v_s,\theta_s,\varphi_s,u_s,x_s)}$ is a Poisson measure on
$[0,1]\times\bb{R}^3\times(0,\pi/2]\times[0,2\pi]\times[0,\infty)\times[0,1]$ with intensity $dsf_s(dv)\beta(\theta)d\theta d\varphi dudx$.
According to Lemma \ref{le}, we know that $f_s(\scr{H}_{w})\ge b$ for all $s\in[0,1]$ and all $w\in \bb{R}^3$.
Then we can get the following lemma.

\blem\label{le2}
For $m\ge1$, we introduce
\[\ca{J}'_m:=\left\{s\in \scr{D}: u_s\le d^\gamma, ~v_s\in \scr{H}_{H_B(V_{s-})},~\theta_s\le K2^{-m},~x_s\le\frac{b}{f_s(\scr{H}_{H_B(V_{s-})})}\right\},\]
where $K=1/(B+c)$ and where $d=a$ (if $\gamma\in(0,1)$) or $d=B+c$ (if $\gamma\in(-1,0]$).
Then we have
\beq\label{rel}
\mathcal{J}_m'\subset\mathcal{J}_m \quad\text{and}\quad \bigcup_{s\in \mathcal {J}_m'}
\Big[s-\Big(\frac{a\theta_s}{4}\Big)^{\delta}, s+\Big(\frac{a\theta_s}{4}\Big)^{\delta}\Big]\subset A_\delta^m.\eeq
\elem

\bpf
We recall that, for all $s\in[0,1]$, $|H_B(V_{s-})|=\abs{\frac{|V_{s-}|\wedge B}{|V_{s-}|}V_{s-}}\le B$ and that $v_s\in \scr{H}_{H_B(V_{s-})}$ implies that $|H_B(V_{s-})-v_s|\ge a$ and $|v_s|\le c$.
Then for all $m\ge1$, for all $s\in\ca{J}_m'$,
we have (recall (\ref{num})) \[|\Delta V_s^B|=\sqrt{\frac{1-\cos\theta_s}{2}}|H_B(V_{s-})-v_s|
\bbd{1}_{\{u_s\le|H_B(V_{s-})-v_s|^\gamma\}}\le\theta_s|H_B(V_{s-})-v_s|\le K2^{-m}(B+c)=2^{-m}.\]
In addition, for all $s\in\ca{J}_m'$, using that $|H_B(V_{s-})-v_s|\geq a$ and
that $1-\cos \theta \geq \theta^2/8$ on $(0,\pi/2]$,
\[|\Delta V_s^B|=\sqrt{\frac{1-\cos\theta_s}{2}}|H_B(V_{s-})-v_s|
\bbd{1}_{\{u_s\le|H_B(V_{s-)}-v_s|^\gamma\}}\ge \frac{a\theta_s}{4}.\]
Indeed, the indicator equals $1$ because we always have $u_s\le d^\gamma\le|H_B(V_{s-})-v_s|^\gamma$
(if $\gamma \in (0,1)$, then $|H_B(V_{s-})-v_s|\geq a$ and $d=a$, while if $\gamma \in (-1,0]$,
then  $|H_B(V_{s-})-v_s|\leq B+c$ and $d=B+c$).
Consequently, for all $m\ge 1$, and all $s \in \ca{J}_m'$, $0<|\Delta V_s^B|\le2^{-m}$: this implies
that $\mathcal{J}_m'\subset\mathcal{J}_m$. Furthermore, for any $\delta>0$,
\[A_\delta^m=\bigcup_{s\in \mathcal {J}_m} [s-|\Delta V_s^B|^{\delta},~s+|\Delta V_s^B|^{\delta}]\supset
\bigcup_{s\in \mathcal {J}_m'} \Big[s-\Big(\frac{a\theta_s}{4}\Big)^{\delta}, ~s+\Big(\frac{a\theta_s}{4}\Big)^{\delta}\Big]\]
as desired.
\epf

\blem\label{le1}
Let $m\ge1$ and $\delta>0$ be fixed. The
random measure \[\mu_m^{\delta}(ds,dy)=\sum_{s\in \mathcal {J}_m'} \dirac_{(s,(a\theta_s/4)^{\delta})}\]
is a Poisson measure on $[0,1]\times(0,\infty)$  with intensity $ds~h_m^{\delta}(y)dy$,
where \[h_m^{\delta}(y)=\frac{8\pi d^\gamma b}{a\delta}\beta\Big(\frac{4}{a}y^{1/\delta}\Big)y^{\frac{1}{\delta}-1}\bbd{1}_{\{y\leq (aK2^{-(m+2)})^{\delta}\wedge(a\pi/8)^\delta\}}.\]
Moreover, we have
\[c_1y^{-1-\frac{\nu}{\delta}}\bbd{1}_{\{y\leq (aK2^{-(m+2)})^{\delta}\wedge(a\pi/8)^\delta\}}\le h_m^{\delta}(y)\leq C_1y^{-1-\frac{\nu}{\delta}}
\bbd{1}_{\{y\leq (aK2^{-(m+2)})^{\delta}\wedge(a\pi/8)^\delta\}},\]
for some constants $0<c_1<C_1$ \emph{(}depending on $B,\delta$\emph{)}.
\elem
\bpf
By Jacod-Shiryaev \cite{jacod1987limit} [Chapter 2, Theorem 1.8], it suffices to check that the compensator of
the random measure $\mu_m^{\delta}(ds,dy)$ is $dsh_m^{\delta}(y)dy$, i.e., for any predictable process $W(s,y)$,
\begin{align*}
&\int_0^t\int_0^{\infty} W(s,y)(\mu_m^{\delta}(ds,dy)-dsh_m^{\delta}(y)dy)\\
=&\int_{0}^{t}\int_{\bb{R}^3}\int_0^{\pi/2}\int_0^{2\pi}\int_0^{\infty}\int_0^1 W(s,(a\theta/4)^{\delta})\times\bbd{1}_{\{v\in \scr{H}_{H_B(V_{s-})},~
\theta\leq K2^{-m},~u\le d^\gamma,~x\le b/f_s(\scr{H}_{H_B(V_{s-})})\}}\\
& \hskip6cm
 \times O(ds,dv,d\theta,d\varphi,du,dx)
- \int_0^t\int_0^{\infty} W(s,y)h_m^{\delta}(y)dsdy
\end{align*}
is a local martingale.
Recalling that $O$ is a Poisson measure with intensity $dsf_s(dv)\beta(\theta)d\theta d\varphi dudx$,
we know that
\begin{align*}
&\int_{0}^{t}\int_{\bb{R}^3}\int_0^{\pi/2}\int_0^{2\pi}\int_0^{\infty}\int_0^1 W(s,(a\theta/4)^{\delta})
\bbd{1}_{\{v\in \scr{H}_{H_B(V_{s-})},~\theta\leq K2^{-m},~u\leq d^\gamma,~x\le b/f_s(\scr{H}_{H_B(V_{s-})})\}}\\
& \hskip6cm \big(O(ds,dv,d\theta,d\varphi,du,dx)-dsf_s(dv)\beta(\theta)d\theta d\varphi dudx\big)
\end{align*}
is a local martingale.
Thus, we only need to prove that
\begin{align*}
 &\int_{0}^{t}\int_{\bb{R}^3}\int_0^{\pi/2}\int_0^{2\pi}\int_0^{\infty}\int_0^1  W(s,(a\theta/4)^{\delta})
\bbd{1}_{\{v\in \scr{H}_{H_B(V_{s-})},~\theta\leq K2^{-m},~u\leq d^\gamma,~x\le b/f_s(\scr{H}_{H_B(V_{s-})})\}}\\
& \hskip10cm
dsf_s(dv)\beta(\theta)d\theta d\varphi dudx\\
 =&\int_0^t\int_0^{\infty} W(s,y)h_m^{\delta}(y)dsdy.
\end{align*}
 Actually,
\begin{align*}
&\int_{0}^{t}\int_{\bb{R}^3}\int_0^{\pi/2}\int_0^{2\pi}\int_0^{\infty}\int_0^1  W(s,(a\theta/4)^{\delta})
\bbd{1}_{\{v\in \scr{H}_{H_B(V_{s-})},~\theta\leq K2^{-m},~u\leq d^\gamma,~x\le b/f_s(\scr{H}_{H_B(V_{s-})})\}} \\
& \hskip10cm
dsf_s(dv)\beta(\theta)d\theta d\varphi dudx\\
=&2\pi d^\gamma b\int_0^t\int_0^{\pi/2} W(s,(a\theta/4)^{\delta})\bbd{1}_{\{\theta\leq K2^{-m}\}}ds\beta(\theta)d\theta.
\end{align*}
Using the substitution $y=(a\theta/4)^\delta$, we conclude that the intensity of $\mu_m^{\delta}$ is indeed $dsh_m^{\delta}(y)dy$.
From (\ref{con}), we can easily get the bounds for $h_m^{\delta}(y)$.
\epf

Now, we give the

\begin{proof}[Proof of Proposition \ref{rfs}.]
We start with (1) and thus fix $\delta\in(0,\nu)$. By Lemma \ref{le1}, we know that the random measure
$\mu_m^{\delta}=\sum_{s\in \ca{J}_m'}\dirac_{(s,(a\theta_s/4)^\delta)}$ is a Poisson measure on $[0,1]\times(0,1)$ with intensity $ds~h_m^{\delta}(y)dy$,
where \[h_m^{\delta}(y)\ge c_1y^{-1-\frac{\nu}{\delta}}\bbd{1}_{\{y\leq (aK2^{-(m+2)})^{\delta}\wedge(a\pi/8)^\delta\}}.\]
Clearly, for all $m\ge1$, $\delta\in(0,\nu)$,
\[\int_0^1\exp\left(2\int_t^1\int_y^1h_m^{\delta}(z)dzdy\right)dt=\infty,\]
since $2\int_t^1(\int_y^1h_m^{\delta}(z)dz)dy \gtrsim 2c_1\frac{\delta^2}{(\nu-\delta)\nu}t^{1-\frac{\nu}{\delta}}$.
Applying Lemma \ref{shep}, we deduce that almost surely, for all $m\ge1$,
\[[0,1] \subset \bigcup_{s\in \mathcal {J}_m'} \Big[s-\Big(\frac{a\theta_s}{4}\Big)^{\delta}, s+\Big(\frac{a\theta_s}{4}\Big)^{\delta}\Big].\]
Consequently, almost surely, \[[0,1]\subset\limsup_{m\rightarrow+\infty}\bigcup_{s\in \mathcal {J}_m'}
\Big[s-\Big(\frac{a\theta_s}{4}\Big)^{\delta}, s+\Big(\frac{a\theta_s}{4}\Big)^{\delta}\Big].\]
Recalling (\ref{def2}) and (\ref{rel}), we deduce that $[0,1]\subset A_\delta$ almost surely.

We next prove (2). We set $m_1=1$. By (1), we have a.s.
$[0,1]\subset A_{\nu(1-1/2)}\subset \bigcup_{m\ge m_1} \widetilde{A}_{\nu(1-1/2)}^m$.
 Hence we can find $m_2>m_1$ such that
\[\mathscr{L}\left(\bigcup_{m_1\le m< m_2}\widetilde{A}_{\nu/2}^m\bigcap[0,1]\right)\ge1-\frac{1}{2}.\]
Similarly, we have almost surely,
$[0,1]\subset A_{\nu(1-1/3)}\subset \bigcup_{m\ge m_2}\widetilde{A}^m_{\nu(1-1/3)}$, therefore we can find $m_3>m_2$
such that
\[\mathscr{L}\left(\bigcup_{m_2\le m< m_3}\widetilde{A}_{\nu(1-1/3)}^m\bigcap[0,1]\right)
\ge1-\frac{1}{2^2}.\]
By induction, we can find an increasing sequence $(m_j)_{j\ge1}$ such that, for all $j\ge2$,
\[\mathscr{L}\left(\bigcup_{m_{j-1}\le m<m_j}\widetilde{A}_{\nu(1-1/j)}^m\bigcap[0,1]\right)\ge1-\frac{1}{2^{j-1}}.\]
So, from the Fatou lemma, we have
\[\mathscr{L}\left(\limsup_{j\rightarrow+\infty}\bigcup_{m_{j-1}\le m< m_j}\widetilde{A}_{\nu(1-1/j)}^m\bigcap[0,1]\right)
\ge\limsup_{j\rightarrow+\infty}\mathscr{L}\left(\bigcup_{m_{j-1}\le m< m_j}\widetilde{A}_{\nu(1-1/j)}^m\bigcap[0,1]\right)\ge1.\]
We now put $\epsilon_m=\frac{1}{j}$ for $m\in[m_{j-1},m_j)$ and note that
\[\limsup_{j\rightarrow+\infty}\bigcup_{m\in[m_{j-1},m_j)}
\widetilde{A}_{\nu(1-\epsilon_m)}^m
=\limsup_{m\rightarrow+\infty}\widetilde{A}_{\nu(1-\epsilon_m)}^m.\]
The conclusion follows.
\epf

\subsection{Study of the H\"{o}lder exponent of $V^B$}
\quad~~We now study the pointwise H\"{o}lder exponent of the localized Boltzmann process $V^B$.
\bdefi\label{def3}
For all $t\in [0,1]$, the index of approximation of $t$ is defined by
\[\delta_t:=\sup\{\delta> 0: t\in A_\delta\}.\]
\edefi
For all $t\in [0,1]$, the index of approximation of $t$ reflects directly the relation between $t$ and jump times of $V^B$.
If $\delta_t$ is large, then $t$ is close to many large jumps of $V^B$.
\brem\label{rem}
Recalling \emph{Remark \ref{rema}} and \emph{Proposition \ref{rfs}}, we see that almost surely, for all $t\in[0,1]$,
$\delta_t=\sup\{\delta>0:t\in A_\delta^*\}$ and $\delta_t\ge\nu$.
\erem
If $t\in\ca{J}$, we know that $h_{V^B}(t)=0$. Then for $t\in[0,1]\setminus\ca{J}$, we claim that the H\"{o}lder exponent is the inverse of the index of approximation.
\bprop\label{he}
Almost surely, for all $t\in[0,1]\setminus\ca{J}$, $h_{V^B}(t)=\frac{1}{\delta_t}$.
\eprop

To prove this claim, we need the following two lemmas. The first lemma, that will give the upper bound
for $h_{V^B}(t)$, can be found in  \cite{jaffard1999multifractal} and is as follows.

\blem\label{l1}
Let $f:\bb{R}\rightarrow \bb{R}^3$ be a function discontinuous on a dense set of points and let $(t_n)_{n\ge 1}$ be a real
sequence converging to some $t$ and such that $f$ has left and right limits at each $t_n$. Then
\[h_f(t)\le \liminf_{n\rightarrow \infty} \frac{\log |f(t_n+)-f(t_n-)|}{\log|t_n-t|}.\]
\elem

For the lower bound of $h_{V^B}(t)$, we will use Lemma \ref{l2} below, that relies on
some ideas of \cite{balancca2014fine}.
We first introduce, for $m>0$,  the following two processes:
\begin{align*}
V_t^{B,m} :=& V_0+\int_{0}^{t}\int_{\bb{R}^3}\int_0^{\pi/2}\int_0^{2\pi}\int_0^{\infty}
{\bf a}(H_B(V_{s-}),v,\theta,\varphi)~
\bbd{1}_{\{u\le |H_B(V_{s-})-v|^\gamma\}}\\
& \hskip6cm
\times\bbd{1}_{\{|{\bf a}(H_B(V_{s-}),v,\theta,\varphi)|\le 2^{-m}\}}N(ds, dv,d\theta,d\varphi,du),\\
Z_t^{B,m} :=& \int_{0}^{t}\int_{\bb{R}^3}\int_0^{\pi/2}\int_0^{2\pi}\int_0^{\infty}
\theta|H_B(V_{s-})-v|
\bbd{1}_{\{u\le |H_B(V_{s-})-v|^\gamma\}}\\
& \hskip6cm
\times\bbd{1}_{\{\frac{\theta}{4}|H_B(V_{s-})-v|\le 2^{-m}\}}N(ds, dv,d\theta,d\varphi,du).
\end{align*}
We can immediately observe that the process $Z_t^{B,m}$ is almost surely increasing as a function of $t$. We also notice that a.s., for all $x,y\in[0,1]$,
\beq\label{com}\abs{V_x^{B,m}-V_y^{B,m}}
\le \abs{Z_x^{B,m}-Z_y^{B,m}}.\eeq
 This comes from the inequality $\frac{\theta}{4}|H_B(V_{s-})-v|\le|{\bf a}(H_B(V_{s-}),v,\theta,\varphi)|\le\theta|H_B(V_{s-})-v|$,
which follows from (\ref{num}).
\blem\label{l2}
There exists some constant $C_{B}>0$, such that
\begin{enumerate}[label=\emph{(\arabic*)}]
\item for all $\delta>\nu$, all $m\ge1$,
\beq\label{eq1}
\bb{P}\left[\sup_{x,y\in[0,1], |x-y|\le 2^{-m}} \abs{V_x^{B,\frac{m}{\delta}}-V_y^{B,\frac{m}{\delta}}}\ge
m2^{-\frac{m}{\delta}}\right]\le C_{B} e^{-m/4},\eeq

\item for all $m\ge1$, all $\lambda\in[0,2^{m}]$,
\beq\label{eq2}
\bb{E}\left[e^{\lambda Z_1^{B,m}}\right]\le e^{C_B\lambda2^{-m(1-\nu)}}.
\eeq
\end{enumerate}
\elem
\bpf
We first prove (\ref{eq1}).
Setting $\lambda=3\times2^{m/\delta}$, recalling (\ref{com}) and that $Z_t^{B,\frac{m}{\delta}}$ is almost surely increasing in $t$, we get
\begin{align*}
 \bb{P}\left[\sup_{x,y\in[0,1], |x-y|\le 2^{-m}}
\abs{V_{x}^{B,\frac{m}{\delta}}-V_{y}^{B,\frac{m}{\delta}}}
\ge m2^{-\frac{m}{\delta}}\right]
&\le\bb{P}\left[\sup_{x,y\in[0,1], |x-y|\le 2^{-m}}
\abs{Z_{x}^{B,\frac{m}{\delta}}-Z_{y}^{B,\frac{m}{\delta}}}
\ge m2^{-\frac{m}{\delta}}\right]\\
&\le\sum_{k=0}^{2^m-1}\bb{P}\left[ \bra{Z_{(k+1)2^{-m}}^{B,\frac{m}{\delta}}-Z_{k2^{-m}}^{B,\frac{m}{\delta}}}
\ge \frac{m2^{-\frac{m}{\delta}}}{3}\right]\\
&\le \sum_{k=0}^{2^m-1} e^{-m} \bb{E}\left[\exp\set{\lambda\bra{Z_{(k+1)2^{-m}}^{B,\frac{m}{\delta}}
-Z_{k2^{-m}}^{B,\frac{m}{\delta}}}}
\right]\\
&=:\sum_{k=0}^{2^m-1} e^{-m} I_k.
\end{align*}
We then set
\[J_k(t):=\bb{E}\left[\exp\set{\lambda\bra{Z_{t+k2^{-m}}^{B,\frac{m}{\delta}}
-Z_{k2^{-m}}^{B,\frac{m}{\delta}}}
}\right].\]
Observe that $I_k=J_k(2^{-m})$. For all $t\geq 0$, we have
\begin{align*}
J_k(t)
= 1+2\pi\bb{E}\Bigg[ & \int_{k2^{-m}}^{t+k2^{-m}}\int_{\bb{R}^3}\int_0^{\pi/2}
\exp\set{\lambda\left(Z_s^{B,\frac{m}{\delta}}-Z_{k2^{-m}}^{B,\frac{m}{\delta}}\right)}
(e^{\lambda\theta|H_B(V_{s})-v|}-1) \\
& \hskip4cm \times |H_B(V_{s})-v|^\gamma\bbd{1}_{\set{\frac{\theta}{4}|H_B(V_{s})-v|\le 2^{-\frac{m}{\delta}}}}
\beta(\theta)d\theta f_s(dv)ds\Bigg].
\end{align*}
From $\lambda\theta\abs{H_B(V_{s})-v}\le4\lambda2^{-m/\delta}=12$, we have $e^{\lambda\theta|H_B(V_{s})-v|}-1 \le C\lambda\theta|H_B(V_{s})-v|$ for some positive constant $C$. Using this estimate and recalling (\ref{con}), we get
\begin{align*}
J_k(t)
\le& 1+C\lambda\bb{E}\Bigg[  \int_{k2^{-m}}^{t+k2^{-m}}\int_{\bb{R}^3}\int_0^{\pi/2}
\exp{\set{\lambda\left(Z_s^{B,\frac{m}{\delta}}-Z_{k2^{-m}}^{B,\frac{m}{\delta}}\right)}}\\
& \hskip6cm \times  \theta^{-\nu} |H_B(V_{s})-v|^{\gamma+1} \bbd{1}_{\set{\frac{\theta}{4}|H_B(V_{s})-v|\le 2^{-\frac{m}{\delta}}}}d\theta
f_s(dv)ds\Bigg].
\end{align*}
Moreover,
\begin{align*}
    |H_B(V_{s})-v|^{\gamma+1}\int_0^{\pi/2}\theta^{-\nu}
    \bbd{1}_{\set{\frac{\theta}{4}|H_B(V_{s})-v|\le 2^{-\frac{m}{\delta}}}}
d\theta\le& C|H_B(V_{s})-v|^{\gamma+1}(|H_B(V_{s})-v|2^{\frac{m}{\delta}})^{\nu-1}\\ \le& C |H_B(V_{s})-v|^{\gamma+\nu} 2^{\frac{m(\nu-1)}{\delta}}.
\end{align*}
Since $\gamma+\nu\in(0,2)$ by assumption, we have $|H_B(V_{s})-v|^{\gamma+\nu}\le C(1+|v|^2+|H_B(V_{s})|^2)$, whence
\begin{align*}
J_k(t)\le 1+C\lambda2^{\frac{m(\nu-1)}{\delta}}\bb{E}\Bigg[ & \int_{k2^{-m}}^{t+k2^{-m}}\int_{\bb{R}^3}
\exp\set{\lambda\left(Z_s^{B,\frac{m}{\delta}}-Z_{k2^{-m}}^{B,\frac{m}{\delta}}\right)}
(1+|H_B(V_{s})|^2+|v|^2) f_s(dv)ds\Bigg].
\end{align*}
Since $|H_B(V_{s})|\le B$, and by conservation of the kinetic energy, we have a.s.
\[\int_{\bb{R}^3}(1+|H_B(V_{s})|^2+|v|^2)f_s(dv)\le1+B^2+m_2(f_0).\]
Using finally that $\lambda~2^{\frac{m(\nu-1)}{\delta}}=3\times 2^{\frac{m\nu}{\delta}}$, we find that
\begin{align*}
J_k(t)\le 1+C_B2^{\frac{m\nu}{\delta}}\bb{E}\Bigg[\int_{k2^{-m}}^{t+k2^{-m}}
\exp{\set{\lambda\left(Z_s^{B,\frac{m}{\delta}}-Z_{k2^{-m}}
^{B,\frac{m}{\delta}}\right)}}ds
\Bigg].
\end{align*}
It follows that $J_k(t)\le1+C_B2^{\frac{m\nu}{\delta}}\int_0^tJ_k(s)ds$. Hence $J_k(t)\le \exp(C_B2^{\frac{m\nu}{\delta}}t)$
by the Gr\"{o}nwall inequality, so that $I_k=J_k(2^{-m}) \le \exp(C_B 2^{-m(1-\frac{\nu}{\delta})}) \le C_{B}$
because $\delta \geq \nu$. Finally,
\[\bb{P}\left[\sup_{x,y\in[0,1], |x-y|\le 2^{-m}}
\abs{V_{x}^{B,\frac{m}{\delta}}-V_{y}^{B,\frac{m}{\delta}}}
\ge m2^{-\frac{m}{\delta}}\right]
\le\sum_{k=0}^{2^m-1} e^{-m} I_k
\le C_{B}~e^{-m}2^m
\le C_{B}~e^{-m/4}.\]
This completes the proof of (\ref{eq1}). We only sketch the proof of (\ref{eq2}), since it is very similar.
First, by It\^{o} Formula,
\begin{align*}
&\bb{E}\left[e^{\lambda Z_t^{B,m}}\right]\\
&=  1+2\pi\bb{E}\Bigg[\int_0^t\int_{\bb{R}^3}\int_0^{\pi/2}e^{\lambda Z_{s}^{B,m}}
\bra{e^{\lambda\theta|H_B(V_{s})-v|}-1}|H_B(V_{s})-v|^\gamma
\bbd{1}_{\set{\frac{\theta}{4}|H_B(V_{s})-v|\le 2^{-m}}} \beta(\theta)d\theta f_s(dv)ds\Bigg].
\end{align*}
Since $\lambda\theta|H_B(V_{s})-v|<4$ (because $\lambda \leq 2^{m}$), a similar computation as previously shows that
\begin{align*}
\bb{E}\left[e^{\lambda Z_{t}^{B,m}}\right]
\le 1+C_B\lambda2^{m(\nu-1)}\bb{E}\Bigg[\int_0^{t}e^{\lambda Z_s^{B,m}}ds\Bigg]
\le 1+C_B\lambda2^{m(\nu-1)}\int_0^t\bb{E}[e^{\lambda Z_{s}^{B,m}}]ds.
\end{align*}
Owing to the Gr\"{o}nwall inequality, we deduce that $\bb{E}[e^{\lambda Z_{t}^{B,m}}]
\leq e^{C_B\lambda2^{m(\nu-1)}t}$. Taking $t=1$, we obtain the conclusion.
\epf

Now, we can proceed to the

\bpf[Proof of Proposition \ref{he}.]
\emph{Upper Bound.} Here we prove that
for all $t\in[0,1]$, it holds that $h_{V^B}(t)\le1/\delta_t$. To this end, we check that for all $\delta>0$, all $t\in A_\delta$, $h_{V^B}(t)\le1/\delta$. Let thus $\delta>0$ and $t\in A_\delta$.
By definition of $A_\delta$, for all $m\ge1$, there exists $t_m\in \mathcal {J}$,
such that $|t_m-t|\le |\Delta V_{t_m}^B|^\delta$ and $|\Delta V_{t_m}^B|\le2^{-m}$. From Lemma \ref{l1},  we directly deduce that
\begin{align*}
h_{V^B}(t)\le \liminf_{m\rightarrow \infty} \frac{\log|\Delta V_{t_m}^B|}{\log|t_m-t|}
\le \liminf_{m\rightarrow \infty} \frac{\log|\Delta V_{t_m}^B|}{\log|\Delta V_{t_m}^B|^{\delta}}
= \frac{1}{\delta}.
\end{align*}
\emph{Lower Bound.}
In this part we show that almost surely, for all $t\in[0,1]\setminus\ca{J}$,
$h_{V^B}(t)\ge1/\delta_t$. To get this, we need to check that for all $\delta>\nu$, if $t\notin A_\delta$,
then $h_{V^B}(t)\ge1/\delta$. Let thus $\delta>\nu$ and $t\notin A_\delta$.

By Lemma \ref{l2}-(1) and Borel-Cantelli's lemma, there almost surely exists $m_0\ge1$
 such that for all $m>m_0$, for all $x,y\in[0,1]$ satisfying $|x-y|\le 2^{-m}$,
\beq\label{cantelli}
|V_x^{B,\frac{m}{\delta}}-V_y^{B,\frac{m}{\delta}}|\le m2^{-\frac{m}{\delta}}.\eeq
Since $t\notin A_\delta$, there exists $m_1>m_0$,
 such that for all $s\in\ca{J}$ satisfying $|\Delta V_s^B| \le 2^{-m_1}$,
we have \beq\label{inequal}|s-t|>|\Delta V_s^B|^{\delta}.\eeq
For all $r\in[0,1]$, we define
\[U_{t,r}^{m_1}:=\sum_{s\in[t\wedge r,t\vee r]\cap\ca{J}}|\Delta V_s^B|~\bbd{1}_{\set{|\Delta V_s^B|>2^{-m_1}}},\] and we observe that \[|V_t^B-V_r^B|\le|V_t^{B,m_1}-V_r^{B,m_1}|+U_{t,r}^{B,m_1}.\]
Since $t\notin\ca{J}$ and since the process $V^B$ has almost surely a finite number of jump greater than $2^{-m_1}$, we can almost surely find $\epsilon_1>0$ such that,
for all $r\in(t-\epsilon_1,t+\epsilon_1)$,
$U_{t,r}^{m_1}=0$.

Next, we put $\epsilon_2=2^{-m_1-1}$. Then for each $r\in (t-\epsilon_2,t+\epsilon_2)$,
we set $m_r=\lfloor \log_2{\frac{1}{|t-r|}}\rfloor>m_1$, for which
$2^{-m_r-1}<|t-r|\le2^{-m_r}$. Then for all $r\in(t-\epsilon_2,t+\epsilon_2)$, we write
\[|V_t^{B,m_1}-V_r^{B,m_1}|\le|V_t^{B,m_r/\delta}-V_r^{B,m_r/\delta}|
+\sum_{s\in[t\wedge r, t\vee r]\cap\ca{J}}|\Delta V_s^B|~\bbd{1}_{\{2^{-\frac{m_r}{\delta}}<|\Delta V_s^B|\le2^{-m_1}\}}.\]
According to (\ref{inequal}), for $s\in [t\wedge r, t\vee r]\cap\ca{J}$, $|\Delta V_s^B|\le2^{-m_1}$ implies that
$|\triangle V_s^B|<|s-t|^{1/\delta}\le|r-t|^{1/\delta}\le2^{-\frac{m_r}{\delta}}$, whence the second term $\sum_{s\in[t\wedge r, t\vee r]\cap\ca{J}}|\Delta V_s^B|~\bbd{1}_{\{2^{-\frac{m_r}{\delta}}<|\Delta V_s^B|\le2^{-m_1}\}}$ vanishes.

To summarize, we have checked that for all $r\in\Big(t-(\epsilon_1\wedge\epsilon_2),~t+(\epsilon_1\wedge\epsilon_2)\Big)$,
\[|V_t^B-V_r^B|\le\Big|V_t^{B,m_r/\delta}-V_r^{B,m_r/\delta}\Big|.\]
Furthermore, since $m_r>m_0$, we conclude from (\ref{cantelli}) that, still for $r\in\Big(t-(\epsilon_1\wedge\epsilon_2),~t+(\epsilon_1\wedge\epsilon_2)\Big)$,
\[|V_t^B-V_r^B|\le m_r2^{-\frac{m_r}{\delta}}\le\frac{2^{1/\delta}}{\log2}
\log\Big(\frac{1}{|t-r|}\Big)|t-r|^{1/\delta}.\]
This implies that $h_{V^B}(t)\ge\frac{1}{\delta}$ and ends the proof.
\epf

\subsection{Hausdorff dimension of the sets $A^*_\delta$}
\quad~~Now, we compute the Hausdorff dimension of $A_\delta^*$, which will be used for giving the spectrum of
singularities and the proof of Proposition \ref{loc} in the next subsection.

\bprop\label{hd} Almost surely, for all $\delta>\nu$,
\[\emph{dim}_H(A_\delta^*)=\frac{\nu}{\delta} ~and~ \ca{H}^{\nu/\delta}(A_\delta^*)=+\infty.\]
\eprop

To check this proposition, we need the \emph{mass transference principle}, proved in \cite{beresnevich2006mass}, Theorem 2 (applied in dimension $k=1$ and
with the function $f(x)=x^\alpha$).

\bprop\label{matp} Let $\alpha\in(0,1)$ be fixed.
Let $\{F_i=[x_i-r_i,~x_i+r_i]\}_{i\in\mathbb{N}}$ be a sequence of intervals in $\bb{R}$ with radius $r_i\rightarrow0$ as $i\rightarrow +\infty$.
Suppose that
\[\mathscr{L}(\limsup_{i\rightarrow +\infty}F_i^\alpha\cap[0,1])=1,\]
where $F_i^\alpha:=[x_i-r_i^\alpha,~x_i+r_i^\alpha]$.
Then,
\[\ca{H}^\alpha(\limsup_{i\rightarrow +\infty}F_i\cap[0,1])=\ca{H}^\alpha([0,1])=+\infty.\]
\eprop

\bpf[Proof of Proposition \ref{hd}.]
%We fix $\delta>\nu$ for the whole proof.
\emph{Lower Bound.} We fix $\delta>\nu$. For all $m\ge1$, we set
\[
N_m:= \sharp\widetilde{\ca{J}}_m=\sharp\{s\in\ca{J}:2^{-m-1}<|\Delta V_s^B|\le2^{-m}\}.\]
We can write $\widetilde{\ca{J}}_m=\{T_1^m,...,T_{N_m}^m\}$, ordered chronologically.
Then we define a sequence $(F_{\delta,j})_{j\ge1}$ of intervals as follows. For $j\ge1$, there is a unique $m\ge1$ and $i\in\{1,2,...,N_m\}$ such that $j=\sum_{k=0}^{m-1}N_k+i$ and write
\[F_{\delta,j}:= \Big[T_i^m-|\Delta V_{T_i^m}^B|^{\delta(1-\epsilon_m)},~ T_i^m+|\Delta V_{T_i^m}^B|^{\delta(1-\epsilon_m)}\Big], \]
where $\epsilon_m$ is defined in Proposition \ref{rfs}.
By this way, we get a sequence of intervals
$(F_{\delta,j})_{j\ge1}$ of radius tending to 0 and such that, for all $\alpha>0$,
$\limsup_{j\rightarrow+\infty}F_{\delta,j}^\alpha=A_{\alpha\delta}^*$
(this is obvious by definition of $A_\delta^*$, see Remark \ref{rema}).
Particularly, taking $\alpha=\frac{\nu}{\delta}\in(0,1)$, we get
\[\limsup_{j\rightarrow+\infty}F_{\delta,j}^{\nu/\delta}=A_\nu^*.\]Thus by Proposition \ref{rfs}-(2),
\[\mathscr{L}\Big(\limsup_{j\rightarrow+\infty}F_{\delta,j}^{\nu/\delta}\cap [0,1]\Big)=1.\]
Consequently, by Proposition \ref{matp}, we have
\[\ca{H}^{\nu/\delta}\Big(\limsup_{j\rightarrow+\infty}F_{\delta,j}\cap[0,1]\Big)=+\infty,\]
that is,
\[\ca{H}^{\nu/\delta}\Big(A_{\delta}^*\cap[0,1]\Big)=+\infty.\]Then $\ca{H}^{\nu/\delta}(A_{\delta}^*)=+\infty$
and $\text{dim}_H(A_\delta^*)\ge \frac{\nu}{\delta}$.

Observing that the family of intervals $F_{\delta,j}^{\nu/\delta}$ does not depend on $\delta$,
we can clearly apply Proposition \ref{matp} simultaneously for all $\delta>\nu$
and we conclude that a.s., for all $\delta>\nu$, $\ca{H}^{\nu/\delta}(A_{\delta}^*)=+\infty$
and $\text{dim}_H(A_\delta^*)\ge \frac{\nu}{\delta}$.

\vskip0.2cm\emph{Upper Bound.} Let $\delta>\nu$ be fixed.
To get the upper bound for $\dim_H(A_\delta^*)$, we show first that a.s.,
$\dim_H(A_\delta)\le\frac{\nu}{\delta}$. For all $m\ge1$,
\[N_m=\sum_{s\in\ca{J}}\bbd{1}_{\{2^{-m-1}<|\Delta V_s^B|\le2^{-m}\}}\le\sum_{s\in\ca{J}} 2^{m+1}|\Delta V_s^B|\bbd{1}_{\{|\Delta V_s^B|\le2^{-m}\}}\le2^{m+1}Z_1^{B,m}.\] This estimate is obtained by using (\ref{com}). Then
\[\bb{P}[N_m\ge m2^{m\nu}]\le\bb{P}[Z_1^{B,m}\ge \frac{1}{2}m2^{m(\nu-1)}].\]
Setting $\lambda=2^{m(1-\nu)}$, we get
\begin{align*}
   \bb{P}[Z_1^{B,m}\ge \frac{1}{2}m2^{m(\nu-1)}]=\bb{P}[\lambda Z_1^{B,m}\ge m/2]\le e^{-\frac{m}{2}}\bb{E}[e^{\lambda Z_1^{B,m}}].
\end{align*}
Since $\lambda=2^{m(1-\nu)} \leq 2^m$, we infer from Lemma \ref{l2}-(2) that
\[\bb{E}[e^{\lambda Z_1^{B,m}}]\le C_B.\]
Hence we obtain
\[\bb{P}[N_m\ge m2^{m\nu}]\le C_B e^{-m/2}.\]
According to the Borel-Cantelli lemma, we know that, almost surely there exists $M>0$ such that,
for all $m>M$, $N_m<m2^{m\nu}$.

Next, by definition of $\widetilde{A}_\delta^k$,
\[\bigcup_{k\ge m}\widetilde{A}_\delta^k\subset\bigcup_{k\ge m}\bigcup_{s\in\widetilde{\ca{J}}_k}[s-2^{-k\delta},s+2^{-k\delta}],\]
so, recalling Definition 1.5, for all $\alpha>0$, and all $m>M$, a.s.,
\begin{align*}
\ca{H}_{2^{-m\delta+1}}^\alpha \Big(\bigcup_{k\ge m}\widetilde{A}_\delta^k\Big)\le 2^\alpha\sum_{k\ge m}N_k2^{-k\delta \alpha}\le 2^\alpha\sum_{k\ge m}k2^{k(\nu-\delta\alpha)}.
\end{align*}
But recalling (\ref{def2}), $A_\delta\subset\bigcup_{k\ge m}\widetilde{A}_\delta^k$, whence,
for all $\alpha>0$, and all $m>M$, a.s.,
\[\ca{H}_{2^{-m\delta+1}}^\alpha (A_\delta)\le 2^\alpha\sum_{k\ge m}k2^{k(\nu-\delta\alpha)}.\]
Consequently,\[\ca{H}^\alpha(A_\delta)=\lim_{m\rightarrow +\infty}
\ca{H}_{2^{-m\delta+1}}^\alpha(A_\delta)\le 2^\alpha\lim_{m\rightarrow +\infty}
\sum_{k\ge m}k2^{k(\nu-\delta\alpha)}.\]
It follows that $\ca{H}^\alpha(A_\delta)=0$ for all $\alpha>\nu/\delta$.
Thus, $\dim_H(A_\delta)\le \nu/\delta$ by Definition \ref{def1}. Since $A_\delta^*\subset A_{\delta^\prime}$
for any $\delta^\prime\in(0,\delta)$, we easily conclude that, a.s.,
\[\dim_H(A_\delta^*)\le\nu/\delta.\]

We have shown that for all $\delta>\nu$, a.s., $\dim_H(A_\delta^*)\le\nu/\delta$.
Using the a.s. monotonicity of $\delta\mapsto A_\delta^*$, it is not hard to conclude that a.s., for all
$\delta>\nu$, $\dim_H(A_\delta^*)\le\nu/\delta$.
\epf

\subsection{Spectrum of singularity of $V^B$}
\quad~~Using Proposition \ref{he}, we can easily get the following relationship between $E_{V^B}(h)$ and $A_\delta^*$.

\bprop\label{enplus}
Almost surely, for all $h>0$, \[
E_{V^B}(h)=\Big(\bigcap_{\delta\in (0,1/h)}A_\delta^*\Big)\setminus\Big(\bigcup_{\delta>1/h}A_\delta^*\Big).
\]
and
\[E_{V^B}(0)=\Big(\bigcap_{\delta\in(0,\infty)}A_\delta^*\Big).\]
\eprop

\brem\label{nostar}
Due to Remark \ref{rema}, Proposition \ref{enplus} also holds when replacing everywhere $A_\delta^*$ by
$A_\delta$.
\erem
We now can finally give the
\bpf[Proof of Proposition \ref{loc}]
We first deal with the case where $h\in(0,1/\nu]$. By Propositions \ref{enplus} and \ref{hd},
\[D_{V^B}(h)=\dim_H\Big(E_{V^B}(h)\Big)\le\dim_H\Big(\bigcap_{\delta\in(0,1/h)} A_\delta^*\Big)
\le\inf_{\delta \in(0,1/h)}\dim_H(A_\delta^*)=h\nu.\]
On the other hand, we observe that (recall that $\delta\mapsto A^*_\delta$ is decreasing)
\[D_{V^B}(h)=\dim_H\Big(E_{V^B}(h)\Big)\ge
\dim_H\Big(A_{{1}/ {h}}^*\setminus(\bigcup_{\delta>1/h}A_\delta^*)\Big).\]
But
\begin{align*}
\ca{H}^{h\nu}\Big(A_{ {1}/ {h}}^*\setminus(\bigcup_{\delta>1/h}A_\delta^*)\Big)
=\ca{H}^{h\nu}(A_{{1}/ {h}}^*)-\ca{H}^{h\nu}\Big(\bigcup_{\delta>1/h}A_\delta^*\Big).
\end{align*}
For all $\delta>{1}/ {h}$, $\dim_H(A_\delta^*)=\frac{\nu}{\delta}<h\nu$, thus $\ca{H}^{h\nu}(A_\delta^*)=0$. Moreover, recalling that $A^*_\delta$ is decreasing when $\delta>\nu$, hence
\[\ca{H}^{h\nu}\Big(\bigcup_{\delta>1/h}A_\delta^*\Big)=0.\]
Next, Proposition \ref{hd} (if $h\nu<1$) and Proposition \ref{rfs} (if $h\nu=1$) imply that
\[\ca{H}^{h\nu}(A_{{1}/ {h}}^*)>0.\]
Consequently, $\dim_H\Big(A_{{1}/ {h}}^*\setminus(\cup_{\delta>1/h}A_\delta^*)\Big)\ge h\nu$, whence finally, $D_{V^B}(h)\ge h\nu$.
We have checked that for $h\in(0,1/\nu]$, it holds that  $D_{V^B}(h)= h\nu$.

When $h=0$, we immediately get, using  Proposition \ref{hd}, that
\[\dim_H\Big(E_{V^B}(0)\Big)=\dim_H\Big(\bigcap_{\delta \in (0,\infty)}A_\delta^*\Big)\le\inf_{\delta \in (0,\infty)}
\frac{\nu}{\delta}=0.\]
Since furthermore $E_{V^B}(0)\supset \cal J$ is a.s. not empty, we conclude that $\dim_H\Big(E_{V^B}(0)\Big)=0$.

Finally, when $h>\frac{1}{\nu}$, we want to show that $\dim_H\Big(E_{V^B}(h)\Big)=-\infty$, i.e. that $E_{V^B}(h)=\emptyset$. This claim immediately follows from Remark \ref{rem} and Proposition \ref{he}, since for all $t\in[0,1]\setminus\ca{J}$, $h_{V^B}(t)=\frac{1}{\delta_t}\le\frac{1}{\nu}$, and for $t\in\ca{J}$, $h_{V^B}(t)=0$.
\epf

\section{Study of the position process}
\quad~~The goal of this last section is to prove Proposition \ref{locpo}.
We thus only consider the case of hard potentials $\gamma \in (0,1)$.
Since $X_t^B=\int_0^tV_s^Bds$, we obviously have a.s., for all
$t \in [0,1]$,
\beq\label{hol}
h_{X^B}(t)\ge 1+h_{V^B}(t).
\eeq
Recall that by Definition, $t\in E^{osc}_{V^B}(h)$ if $h_{X^B}(t)> 1+h_{V^B}(t)$
and $t \in  E^{cusp}_{V^B}(h)$ if $h_{X^B}(t)=1+h_{V^B}(t)$.
Inspired by the ideas of Balan\c{c}a \cite{balancca2014fine}, we will  prove several technical lemmas to get
Proposition \ref{locpo}.

\subsection{Preliminaries}
For any $m>0$ and any interval $[r,t]\subset [0,1]$, we set
\beq\label{set2}
H_{[r,t]}^{m}:=\sharp\{s\in [r,t]\cap\ca{J}:|\Delta V_s^B|\ge 2^{-m}\}.
\eeq

\blem\label{le3}
For any $m\ge1$ and any interval $[r,t]\subset [0,1]$,
\begin{enumerate}[label=\emph{(\arabic*)}]
  \item we have
  \[H_{[r,t]}^{m}\le R_{[r,t]}^{m},\]
  where $R_{[r,t]}^{m}=\int_r^t\int_{\bb{R}^3}
    \int_0^{\pi/2}\int_0^{2\pi}\int_0^{\infty}\bbd{1}_{\{\theta(B+|v|)
    \ge2^{-m}\}}\bbd{1}_{\{u\le(B+|v|)^\gamma\}}
    N(ds,dv,d\theta,d\varphi,du)$;
  \item and, with $a>0$ introduced in Lemma \ref{le} (this actually holds true for any value of $a>0$),
      \[H_{[r,t]}^{m}\ge S_{[r,t]}^{m},\] where $S_{[r,t]}^{m}=\int_r^t\int_{\bb{R}^3}\int_0^{\pi/2}\int_0^{2\pi}\int_0^\infty
    \bbd{1}_{\{|v-H_B(V_{s-})|\ge a\}}\bbd{1}_{\{\theta\ge2^{-m+2}/a\}}\bbd{1}_{\{u\le a^\gamma\}}N(ds,dv,d\theta,d\varphi,du)$.
    \end{enumerate}
\elem
\bpf
By definition of $V^B$, see (\ref{newprocess}), we have
\[
H_{[r,t]}^{m} = \int_r^t\int_{\bb{R}^3}\int_0^{\pi/2}\int_0^{2\pi}\int_0^{\infty}
\bbd{1}_{\{|{\bf a}(H_B(V_{s-}),v,\theta,\varphi))| \geq 2^{-m}\}}\bbd{1}_{\{u\le |H_B(V_{s-})-v|^\gamma\}} N(ds,dv,d\theta,d\varphi,du).
\]
Then the claims immediately follow from
$\frac{\theta}{4}\abs{H_B(V)-v}\le\abs{{\bf a}\big(H_B(V),v,\theta,\varphi\big)}\le\theta(B+|v|)$,
see (\ref{num}),
and $\abs{H_B(V)-v}^\gamma\le(B+|v|)^\gamma$.
\epf
\brem\label{rem1}
Glancing at their definitions, it is clear that $S_{[r,t]}^{m}$ and $R_{[r,t]}^{m}$ are $\scr{F}_t$-measurable,
that $R_{[r,t]}^{m}$ is independent of $\scr{F}_r$ and is a Poisson variable with parameter $(\lambda_{[r,t]}^{m})$, where
\begin{equation}
\label{eq11}
\lambda_{[r,t]}^{m}=\int_r^t\int_{\bb{R}^3}\int_0^{\pi/2}\int_0^{2\pi}\int_0^{\infty}
    \bbd{1}_{\{\theta(B+|v|)\ge2^{-m}\}}
    \bbd{1}_{\{u\le(B+|v|)^\gamma\}}dsf_s(dv)\beta(\theta)d\theta d\varphi du.
    \end{equation}
Using (\ref{con}) and that $m_2(f_s)=m_2(f_0)$ for all $s\in [0,1]$, one easily checks that
there exists a constant $C_B>0$ such that $\lambda_{[r,t]}^{m}\le C_B2^{m\nu }|t-r|$ for all $m>0$
and all $0\leq r \leq t \leq 1$ .
\erem

\subsection{Refined study of the jumps}

\quad~~The goal of this part is to prove the following crucial fact.

\blem\label{jum}
Fix $\epsilon>0$ and set $\alpha=\nu(1-2\epsilon)$ and $\beta=\nu(1+4\epsilon)$.
Almost surely, there exists $M\ge1$, such that for all $m\ge M$, for all $t\in[0,1]$, there exists
$t_m\in B(t,2^{-m\alpha})$ such that $|\Delta V_{t_m}^B|\ge2^{-m}$ and there is no other jump of size greater than
$2^{-m(1+\epsilon)}$ in $B(t_m, 2^{-m\beta}/3)$.
\elem

We start with an intermediate result.

\blem\label{prob}
Fix $\epsilon>0$, $\alpha=\nu(1-2\epsilon)$ and $\beta=\nu(1+4\epsilon)$. For any interval $I=[t_0,t_3)\subset [0,1]$ with length $2^{-m\beta}$,  divide $I=[t_0,t_1)\cup[t_1,t_2)\cup[t_2,t_3)$ into three consecutive intervals with length $2^{-m\beta}/3$. Consider the event \[
A_I^{m,\epsilon}=\{H_{[t_0,t_1)}^{m(1+\epsilon)}=0\}
\cap\{H_{[t_1,t_2)}^{m(1+\epsilon)}=H_{[t_1,t_2)}^m=1\}\cap
\{H_{[t_2,t_3)}^{m(1+\epsilon)}=0\}.
\]
There exist some constants
$c_B>0$ and $m_\e>0$ such that, for all $m\geq m_\e$, all intervals $I\subset[0,1]$ with length $2^{-m\beta}$,
\beq\label{probability}
\bb{P}[A_I^{m,\epsilon}|\scr{F}_{t_0}]\ge c_B2^{-4m\nu\epsilon}.
\eeq
\elem

\bpf
We introduce $A_1=\{H_{[t_0,t_1)}^{m(1+\epsilon)}=0\}$, $A_2=\{H_{[t_1,t_2)}^{m(1+\epsilon)}=H_{[t_1,t_2)}^m=1\}$
and $A_3=\{H_{[t_2,t_3)}^{m(1+\epsilon)}=0\}$, so that $A_I^{m,\e}=A_1\cap A_2\cap A_3$.

\vskip 0.2cm\emph{Step 1.}
First we write, since  $A_1\cap A_2 \in \scr{F}_{t_2}$,
\[\bb{P}[A_{I}^{m,\epsilon}|\scr{F}_{t_0}]
 =\bb{E}\left[\bbd{1}_{A_1\cap A_2}\bb{P}[A_3\big|\scr{F}_{t_2}]\big|\scr{F}_{t_0}\right].\]
But using Lemma \ref{le3} and Remark \ref{rem1},
\begin{align*}
    \bb{P}[A_3\big|\scr{F}_{t_2}] =&\bb{P}\Big[H_{[t_2,t_3)}^{m(1+\epsilon)}=0\Big|\scr{F}_{t_2}\Big]
    \ge\bb{P}\Big[R_{[t_2,t_3)}^{m(1+\epsilon)}=0\Big|\scr{F}_{t_2}\Big]
    =\exp(-\lambda_{[t_2,t_3)}^{m(1+\epsilon)}) \geq \frac12
\end{align*}
for all $m$ large enough (depending only on $\e$), since $\lambda_{[t_2,t_3)}^{m(1+\epsilon)} \leq C_B 2^{m\nu}2^{-m\beta}/3
\leq C_B 2^{-3m\e}/3$.
Consequently, for all $m$ large enough (depending only on $\epsilon>0$), we a.s. have
\beq\label{estimation1}
 \bb{P}[A_{I}^{m,\epsilon}|\scr{F}_{t_0}]
 \ge\frac{1}{2}\bb{P}[A_1\cap A_2|\scr{F}_{t_0}].
\eeq
\vskip 0.2cm\emph{Step 2.}
We next write
\[\bb{P}[A_1\cap A_2|\scr{F}_{t_0}]
=\bb{E}\Big[\bbd{1}_{A_1}\bb{P}[
A_2|\scr{F}_{t_1}]\Big|\scr{F}_{t_0}\Big].\]
But using again Lemma \ref{le3},
\begin{align*}
A_2=&\{H_{[t_1,t_2)}^m\ge1\} \setminus\{H_{[t_1,t_2)}^{m(1+\epsilon)}\ge2\}  \supset
\{S_{[t_1,t_2)}^m\ge1\}\setminus\{R_{[t_1,t_2)}^{m(1+\epsilon)}\ge2\}.
\end{align*}
Thus,
\[
    \bb{P}\big[A_2|\scr{F}_{t_1}\big]
   \ge\bb{P}\big[S_{[t_1,t_2)}^m\ge1\big|\scr{F}_{t_1}\big]- \bb{P}\big[R_{[t_1,t_2)}^{m(1+\epsilon)}\ge2\big|\scr{F}_{t_1}\big].
    \]
First, by Remark \ref{rem1},
\begin{align*}
\bb{P}\big[R_{[t_1,t_2)}^{m(1+\epsilon)}\ge2\big|\scr{F}_{t_1}\big]
    =&1-\Big(1+\lambda_{[t_1,t_2)}^{m(1+\epsilon)}\Big)
    \exp\Big(-\lambda_{[t_1,t_2)}^{m(1+\epsilon)}\Big)
    \le \Big(\lambda_{[t_1,t_2)}^{m(1+\epsilon)}\Big)^2\le C_B2^{-6m\nu\epsilon}.
\end{align*}
Next, we put $Y_t:=S_{[t_1,t)}^m$ for $t\ge t_1$ and observe, according to It\^{o}'s Formula, that
\begin{align*}
\bbd{1}_{\{Y_t=0\}}
&=1+\int_{t_1}^t\int_{\bb{R}^3}\int_0^{\pi/2}\int_0^{2\pi}\int_0^\infty
    \bbd{1}_{\{|v-H_B(V_{s-})|\ge a\}}\bbd{1}_{\{u\le a^\gamma\}}\bbd{1}_{\{\theta\ge2^{-m+2}/a\}}\\
   &\hskip6cm \times\big(\bbd{1}_{\{Y_{s-}+\triangle Y_s=0\}}-\bbd{1}_{\{Y_{s-}=0\}}\big)N(ds,dv,d\theta,d\varphi,du)\\
&=1-\int_{t_1}^t\int_{\bb{R}^3}\int_{2^{-m+2}/a}^{\pi/2}\int_0^{2\pi}\int_0^{a^\gamma}
   \bbd{1}_{\{|v-H_B(V_{s-})|\ge a\}}\bbd{1}_{\{Y_{s-}=0\}}N(ds,dv,d\theta,d\varphi,du).
\end{align*}
Hence, for all $t\geq t_1$,
\begin{align*}
\frac{d}{dt}\bb{E}\big[\bbd{1}_{\{Y_t=0\}}\big|\scr{F}_{t_1}\big]
=-\bb{E}\left[\int_{\bb{R}^3}\int_{2^{-m+2}/a}^{\pi/2}\int_0^{2\pi}\int_0^{a^\gamma}
   \bbd{1}_{\{|v-H_B(V_{t})|\ge a\}}\bbd{1}_{\{Y_t=0\}}f_t(dv)\beta(\theta)d\theta d\varphi du
   \Big|\scr{F}_{t_1}\right].
\end{align*}
Using (\ref{con}) and Lemma \ref{le} (which implies that $f_s(\{v\in\bb{R}^3:\abs{v-H_B(V_{s})}\ge a\})\ge b>0$
a.s. for all $s \in [0,1]$), we easily deduce see that
\[\frac{d}{dt}\bb{E}\big[\bbd{1}_{\{Y_t=0\}}\big|\scr{F}_{t_1}\big]\le -\kappa2^{m\nu}\bb{E}\big[\bbd{1}_{\{Y_{t}=0\}}
\big|\scr{F}_{t_1}\big],\]
for some positive constant $\kappa$.
Integrating this inequality, we deduce that a.s., for all $t\geq t_1$,
\[\bb{E}\big(\bbd{1}_{\{Y_t=0\}}\big|\scr{F}_{t_1}\big)\le \exp\{-\kappa2^{m\nu}(t-t_1)\}.\]
Consequently,
\[
\bb{P}\big[S_{[t_1,t_2)}^m\ge1\big|\scr{F}_{t_1}\big]
    =1-\bb{E}\big(\bbd{1}_{\{Y_{t_2}=0\}}\big|\scr{F}_{t_1}\big)
\ge1-\exp\{-\kappa2^{m\nu}(t_2-t_1)\}
=1-\exp\{-\kappa2^{-4m\nu\epsilon}/3\}.
\]
Finally, for all $m$ large enough (depending only on $\e$), we a.s. have
\[\bb{P}\big[A_2|\scr{F}_{t_1}\big]
   \ge 1-\exp\{-\kappa2^{-4m\nu\epsilon}/3\}
   -C_B2^{-6m\nu\epsilon}\ge c_B2^{-4m\nu\epsilon}.\]
\vskip 0.2cm\emph{Step 3.}
Finally, exactly as Step 1, we obtain that for all $m$ large enough,
$$\bb{P}[A_1|\scr{F}_{t_0}]\ge\frac{1}{2}.$$
\vskip 0.2cm\emph{Step 4.} It suffices to
gather Steps 1, 2 and 3 to conclude the proof.
\epf

\bpf [Proof of Lemma \ref{jum}]
We thus fix $\epsilon>0$ and consider $\alpha$ and $\beta$ as in the statement.
For $m>0$, we introduce the notation $r_m=2^{-m\beta}/3$.
We also introduce the number $q_m^2:=\lfloor2^{m(\beta-\alpha)}\rfloor$, the length
$\ell_m:= q_m^2 2^{-m\beta}$ (we have
$\ell_m \leq 2^{-m\alpha}$ and $\ell_m \simeq  2^{-m\alpha}$) and the number $q_m^1:=\lfloor 1/ \ell_m\rfloor +1$
(we have $q_m^1 \simeq 2^{m\alpha}$).
We consider a covering of $[0,1]$ by $q_m^1$ consecutive intervals
$I_1^m,\dots,I_{q_m^1}^m$ with length $\ell_m$.
Next, we divide each $I_i^m$ into $q_m^2$ consecutive intervals
$I_{i,1}^m, \dots, I_{i,q_m^2}^m$ with length  $2^{-m\beta}$.
Finally, we divide each $I_{i,j}^m$ into three consecutive intervals with length $r_m$, writing
$I_{i,j}^m=[t_{i,j}^m,t_{i,j}^m+r_m)\cup[t_{i,j}^m+r_m,t_{i,j}^m+2r_m)\cup [t_{i,j}^m+2r_m,t_{i,j+1}^m)$.
We consider the event
\begin{align*}
A_{i,j}^m=&\{H_{[t_{i,j}^m,t_{i,j}^m+r_m)}^{m(1+\epsilon)}=0\}\cap
\{H_{[t_{i,j}^m+r_m,t_{i,j}^m+2r_m)}^{m(1+\epsilon)}
=H_{[t_{i,j}^m+r_m,t_{i,j}^m+2r_m)}^m=1\}
\cap
\{H_{[t_{i,j}^m+2r_m,t_{i,j+1}^m)}^{m(1+\epsilon)}=0\}.
\end{align*}
According to Lemma \ref{prob}, we know that if $m$ is large enough (depending only on $\e$), a.s.,
for all $i,j$
\begin{equation}\label{jab}
\bb{P}[A_{i,j}^m|\scr{F}_{t_{i,j}^m}]\ge c_B2^{-4m\nu\epsilon}.
\end{equation}
We now consider, for each $i$, the event
\[K_{m,i}=\bigcap_{j=1}^{q_m^2}(A_{i,j}^m)^c.\]
Then, we easily deduce from (\ref{jab}), together with the fact that
$A_{i,1}^m,\dots,A_{i,j-1}^m \in \scr{F}_{t_{i,j}^m}$ for all $j=1,\dots,q_m^2-1$, that
\begin{align*}
 \bb{P}(K_{m,i}) \leq (1-c_B 2^{-4m\nu\epsilon})^{q_m^2} \leq  (1-c_B 2^{-4m\nu\epsilon})^{2^{m(\beta-\alpha)}-1}.
\end{align*}
Thus for $m$ large enough (depending only on $\e$), we conclude that
\begin{align*}
 \bb{P}(K_{m,i}) \leq \exp \big(-c_B2^{-4m\nu\epsilon}2^{m(\beta-\alpha)}\big)
= \exp \big(-c_B 2^{2 m\nu\epsilon}\big).
\end{align*}
Next, we introduce the event $K_m=\bigcup_{i=1}^{q_m^1}K_{m,i}$. Clearly, for $m$ large enough,
(allowing the value of the constant $c_B>0$ to change)
\[\bb{P}(K_m)\le q_m^1 \exp(- c_B 2^{2m\nu\epsilon}) \leq  \exp(- c_B 2^{2m\nu\epsilon}).
\]
Finally, using the Borel-Cantelli lemma, we conclude that there a.s. exists $M>0$
such that for all $m\geq M$, the event $K_m^c$ is realized
(whence for all $i=1,\dots,q_m^1$, there is $j\in \{1,\dots,q_m^2\}$ such that $A^m_{i,j}$ is realized).
This implies that a.s., for all $m\geq M$, for all $t\in [0,1]$, considering $i \in \{1,\dots,q_m^1\}$ such that
$t \in I_i^m$ and $j\in \{1,\dots,q_m^2\}$ such that $A^m_{i,j}$ is realized,
$V^B$ has exactly one jump greater than $2^{-m(1+\e)}$ in the time interval $I_{i,j}^m$,
this jump is greater than $2^{-m}$ and happens at some time $t_m$ located in the middle of
$I_{i,j}^m$ (more precisely, the distance between $t_m$ and the extremities of $I_{i,j}^m$
is at least $r_m$). We clearly have $|t_m-t|\leq \ell_m \leq 2^{-m\alpha}$,
$|\Delta V^B_{t_m}| \geq 2^{-m}$, and $V^B$ has no other jump of size greater than
$2^{-m(1+\e)}$ in $B(t_m,r_m) \subset I_{i,j}^m$. The proof is complete.
\epf

\subsection{Uniform bound for the H\"older exponent of $X^B$}

\quad~~We show here that $D_{X^B}(h)=-\infty$ for all $h>1+1/\nu$.
We  use a general result for primitives of discontinuous functions.
It based on Proposition 1 in \cite{arneodo1998singularity}, recalled in the following lemma.

\blem\label{wavelet}
Let $\eta >0$ and let $N> \eta$ be an integer. Let $g:\bb{R}\rightarrow\bb{R}$ be a locally bounded function
and let $\psi$ be a $C^\infty$ compactly supported function
with its $N$ first moments vanishing, i.e. $\int_{\bb{R}} x^k\psi(x)dx=0$ for $k=0,\dots,N-1$.
The wavelet transform of $g$ is defined by
\beq\label{wa}  W_\psi(g,a,b) = \frac{1}{a}\int_{\bb{R}} g(t) \; \psi \left( \frac{t-b}{a}\right) dt . \eeq
If $g\in C^\eta (t_0)$, then there exists a constant $C>0$ such that for all $a>0$, all $b\in [t_0-1,t_0+1]$,
\beq \label{we}   \qquad |W_\psi(g,a,b) | \le C   \left(  a^\eta + \left| {t_0-b}\right|^\eta  \right).\eeq
\elem

\newcommand{\osc}{\mbox{osc}}

Now, we give the following general result. For any function $g:\mathbb{R}\to \mathbb{R}$, and any interval $I\subset \mathbb{R}$, we set
$$\osc_I(g) = \sup_{x\in I} g(x) - \inf_{x\in I}g(x).$$
\blem\label{upbound}
Let $g:[0,\infty)\rightarrow\bb{R}$ be a c\`{a}dl\`{a}g function, discontinuous on a dense set of points,
let $G(t)=\int_0^t g(s)ds$. Let $t>0$ and let $(t_m)_{m\ge1}$ be a sequence of discontinuities of the
function $g$ converging to $t$.
For all $s\in\bb{R}$, all $m\ge1$, we define
\beq\label{rwf}
g_m(s)=g(s) - J_m\bbd{1}_{\{s\ge t_m\}},\eeq
where $J_m=g(t_m+)-g(t_m-)$.
Assume that for all $m\geq 1$, there exist $r_m>0$ and $\delta_m>0$ such that
\beq\label{rrr}
\osc_{[t_m-r_m,~t_m+r_m]}(g_m) \le\delta_m  \
\emph{and} \ \lim_{m\rightarrow+\infty}\frac{\delta_m}{|J_m|}=0.
\eeq
Then
\beq\label{hp}
h_G(t)\le\liminf_{m\rightarrow+\infty}
\frac{\log\Big(r_m|J_m|\Big)}{\log\Big(|t_m-t|+r_m\Big)}.
\eeq
\elem

\bpf
Let $\varphi$ be a positive   $C^{\infty}$ function, supported on $[0,1]$ satisfying
$\int_\mathbb{R} \varphi(x) dx=1$.

For $k\geq 1$, let $\psi_k(t)=\varphi^{(k)}(t)$, it is clear that $\psi_k$ is $C^{\infty}$, supported on $[0,1]$
and that its $k$ first moments vanish, so it is a {\it wavelet}.

We now pick an integer $N$ such that $N-2$ is larger than the right hand side of (\ref{hp}), and we denote
 by $c_N(a,b) := W_{\psi_N}(g,a,b) $ and $C_{N+1}(a,b):=W_{\psi_{N+1}}(G,a,b) $ the wavelet transforms of $g$ and $G$ using the wavelet $\psi_N$ and $\psi_{N+1}$, respectively.
An integration by parts shows that
\beq\label{intpart} c_N(a,b) = -\frac{1}{a} C_{N+1} (a,b) . \eeq
We fix $\theta \in (0,1)$ such that $\psi_{N-1}(\theta) > 0$.
It follows from (\ref{rwf}) that $   c_N(r_m,t_m-\theta r_m)= P_m +Q_m$, where
\begin{align*}
P_m=& \frac{1}{r_m}\int_{-\infty}^{+\infty} J_m\bbd{1}_{\{s\ge t_m\}}
\psi_N \left( \frac{s-t_m+\theta r_m}{r_m}\right) ds  = \frac{J_m}{r_m}\int_{t_m}^{+\infty}\psi_N \left( \frac{s-t_m+\theta r_m}{r_m}\right) ds  =-J_m\psi_{N-1} (\theta)
\end{align*}
and
$$Q_m=\frac{1}{r_m}\int_{-\infty}^{+\infty} g_m(s) \psi_N \left( \frac{s-t_m+\theta r_m}{r_m}\right)
ds=\frac{1}{r_m}\int_{-\infty}^{+\infty} (g_m(s)-g_m(t_m)) \psi_N \left( \frac{s-t_m+\theta r_m}{r_m}\right)
ds,$$
where we used that    $\psi_N$ has a vanishing integral. Observing that
\[supp \; \left( \psi_N \left( \frac{\cdot  -t_m+\theta r_m}{r_m}\right)\right)  \subset [t_m -r_m, t_m+ r_m ]\]
and recalling (\ref{rrr}), we deduce that $|Q_m| \leq 2\|\psi_N\|_{\infty} \delta_m$.
As a conclusion,
$$
|c_N(r_m,t_m-\theta r_m)| \geq |P_m|- |Q_m| \geq \psi_{N-1} (\theta) |J_m| - 2\|\psi_N\|_{\infty}  \delta_m \geq c |J_m|
$$
for all $m$ large enough, since $\lim_{m\rightarrow+\infty}\frac{\delta_m}{|J_m|}=0$ by assumption.
Then we obtain according to (\ref{intpart}),
\beq\label{wavcoef}|C_{N+1}(r_m,t_m-\theta r_m)|\ge c r_m |J_m|.\eeq

Assume that  $G\in C^\eta (t)$ for some $\eta >\liminf_{m\rightarrow+\infty}
[\log (r_m|J_m|)/[\log(|t_m-t|+r_m)]$. We   apply Lemma \ref{wavelet} with $g = G$, $\psi = \psi_{N+1}$, $a= r_m$, $b=t_m-\theta r_m$. Hence, there is a constant $C$ such that for all $m$,
\[|C_{N+1} (r_m ,t_m-\theta r_m)|  \leq  C \; (r_m^\eta + | t -t_m+\theta r_m |^\eta) \leq C (r_m+ | t -t_m|)^\eta
. \]
This contradicts (\ref{wavcoef}), so necessarily \eqref{hp} hold true.
\epf

We next apply this lemma to our position process to get a uniform upper bound for all pointwise H\"older exponents of $X^B$.

\bprop\label{uniform}
Almost surely, for all $t\in[0,1]$, the H\"older exponent of $X^B$ satisfies
\beq
h_{X^B}(t)\le 1+\frac{1}{\nu}.
\eeq
\eprop

\bpf
We fix $\epsilon>0$ and set $\alpha=\nu(1-2\epsilon)$ and $\beta=\nu(1+4\epsilon)$.
We  show that a.s., $h_{X^B}(t)\le (1+\beta)/\alpha$ for all $t\in [0,1]$. This clearly suffices
since $\e>0$ can be chosen arbitrarily small.

Lemma \ref{jum} shows that there a.s. exists $M>0$, such that for all $m\ge M$, for all $t\in[0,1]$, there exists
$t_m\in B(t,2^{-m\alpha})$ such that $|\Delta V_{t_m}^B|\ge2^{-m}$ and such that
there is no other jump of size greater than
$2^{-m(1+\epsilon)}$ in $B(t_m, r_m)$, with $r_m:=2^{-m\beta}/3$.

We now fix $t \in [0,1]$.  Up to extraction, one can assume that   the first coordinate  $\widetilde V_s^B$ of the three-dimensional vector $V_s^B$ satisfies $|\Delta \widetilde V_{t_m}^{B}| \ge 2^{-m}/3$.
We now apply Lemma \ref{upbound} with
$g=\widetilde V^{B}$ and $r_m=2^{-m \beta}/3$. We thus introduce $g_m(s)=g(s)-\Delta \widetilde V_{t_m}^{B}\bbd{1}_{\{s\ge t_m\}}$.
Since $V^B$ (and so $\widetilde V_s^B$) has no jump with size greater than $2^{-m(1+\epsilon)}$ within the interval  $B(t_n,r_n)=(t_m-r_m,t_m+r_m)$, we observe that
\[\osc_{B(t_n,r_n)} (g_m) \leq 2\times
\sup_{x,y\in[0,1], |x-y|\le 2^{-m\beta}}|V_x^{B,m(1+\e)} - V_y^{B,m(1+\e)}|.\]
Next, using  Lemma \ref{l2}-(1) (with $\delta=\beta/(1+\e)>\nu$) and the Borel-Cantelli Lemma, we deduce that
there is a.s. $M'>0$ such that, for all $m\geq M'$, all $0<x<y<1$ with $|x-y|<2^{-m\beta}$,
$|V^{B,m(1+\e)}_x-V^{B,m(1+\e)}_y| \leq m \beta 2^{-m(1+\e)}$. That is,
\[ \osc_{B(t_n,r_n)} (g_m) \leq 2m \beta 2^{-m(1+\e)}.\]
Since furthermore $\lim_{m\rightarrow+\infty}\frac{2m \beta 2^{-m(1+\e)}}{|\Delta \widetilde V_{t_m}^{B}|}\le
\lim_{m\rightarrow+\infty}\frac{2m \beta 2^{-m(1+\e)}}{2^{-m}/3}=0$, we can apply Lemma \ref{upbound}
with $\delta_m=2m \beta 2^{-m(1+\e)}$:
\[h_{X^B}(t)\le\liminf_{m\rightarrow+\infty}
\frac{\log\Big(r_m|\Delta \widetilde V_{t_m}^{B}|\Big)}{\log(|t_m-t|+r_m)}\le\liminf_{m\rightarrow+\infty}
\frac{\log\Big(2^{-m(1+\beta)}/9\Big)}{\log (2.2^{-m\alpha})}=\frac{1+\beta}{\alpha}.\]
We used that $r_m|\Delta \widetilde V_{t_m}^{B}| \ge (2^{-m}/3)(2^{-m\beta}/3)$ and that
$|t_m-t| + r_m \le 2^{-m\alpha}+2^{-m\beta}/3\leq 2.2^{-m\alpha}$.
This ends the proof.
\epf

\subsection{Study of the oscillating singularities of $X^B$}

\quad~~To characterize more precisely the set of oscillating times,
we first give the following lemma.
\blem\label{jump}
Let $\delta>\nu$, $\epsilon>0$ and $k\in\bb{N}$ satisfy $\delta>\nu(1+\epsilon)(k+1)/k$.
For all $m\in\bb{N}$, let $(I_{j}^m)_{j=1,\dots,\lfloor 2^{m\delta}\rfloor+1}$ be the
covering of $[0,1]$ composed of successive intervals of
length $2^{-m\delta}$. Almost surely, there exists $M\ge 1$ such that for all $m\ge M$, for all
$j=1,\dots, \lfloor 2^{m\delta} \rfloor$, recalling (\ref{set2}),
\beq\label{iq}
H_{I_{j}^m\cup I^m_{j+1}}^{m(1+\epsilon)} \le k,
\eeq
\elem
\bpf
Using Lemma \ref{le3} and Remark \ref{rem1},
\[\bb{P}\Big(H_{I_{j}^m\cup I_{j+1}^m}^{m(1+\epsilon)}>k\Big)\le\bb{P}\Big(R_{I_{j}^m\cup I_{j+1}^m}^{m(1+\epsilon)}>k\Big)
\le \sum_{\ell=k+1}^{+\infty}\frac{(\lambda_{I_{j}^m\cup I_{j+1}^m}^{m(1+\epsilon)})^\ell}{\ell!}
e^{-\lambda_{I_{j}^m\cup I_{j+1}^m}^{m(1+\epsilon)}} \le  (\lambda_{I_{j}^m\cup I_{j+1}^m}^{m(1+\epsilon)})^{k+1},\]
where the value of $  \lambda_{I_{j}^m\cup I_{j+1}^m}^{m(1+\epsilon)}$ is given by equation \eqref{eq11}.
But, since the length of $I_{j}^m\cup I_{j+1}^m$ is $2.2^{-m\delta}$,
we apply the upper bound found for $\lambda^m_{[r,s]}$ in Remark \ref{rem1} to get $\lambda_{I_{j}^m\cup I_{j+1}^m}^{m(1+\epsilon)}\leq 2 C_B 2^{m\nu(1+\e)-m\delta}$,
so that
\[\bb{P}\Big(H_{I_{j}^m\cup I_{j+1}^m}^{m(1+\epsilon)}>k\Big)\le 2C_B 2^{1+m(k+1)(\nu(1+\epsilon)-\delta)}.\]
Consequently,
\[\bb{P}\Big(\bigcup_{j=1}^{\lfloor 2^{m\delta}\rfloor+1} \Big\{H_{I_{j}^m\cup I_{j+1}^m}^{m(1+\epsilon)}>k\Big\}\Big)
\le 2C_B 2^{m\delta}2^{m(k+1)(\nu(1+\epsilon)-\delta)}= 2C_B 2^{-m k(\delta-\nu(1+\e)(k+1)/k)}.\]
By assumption, this is the general term of a convergent series.
We conclude thanks to the Borel-Cantelli lemma.
\epf

We first study the case where  $h \in [0,1/(2\nu))$.

\blem\label{rela}
Almost surely, for all $h\in[0,1/(2\nu))$, $E^{osc}_{V^B}(h)=\emptyset$.
\elem

\bpf
According to (\ref{hol}), it is sufficient to check that for $h\in[0,1/(2\nu)]$, for all $t\in E_{V^B}(h)$,
$h_{X^B}(t)\le 1+h$.
We fix $\e>0$  so small   that there exists $\delta \in (\max\{2\nu(1+\e),1/(h+\e)\}, 1/h)$.
Next, we  fix $t \in E_{V^B}(h)$. By Remark \ref{nostar}, we know that $t \in A_{1/(h+\e)}$.
Hence for all $n\geq 1$, we can find $m_n\geq n$ and $t_n \in \widetilde {\cal J}_{m_n}$ (that is
$|\Delta V^B_{t_n}| \in (2^{-m_n-1},2^{-m_n}]$) such that $|t_n-t| \leq |\Delta V^B_{t_n}|^{1/(h+\e)}\leq
2^{-m_n/(h+\e)}$.
Applying Lemma \ref{jump} with $k=1$ (since $\delta>2\nu(1+\e)$),
we deduce that $V^B$ has no other jump of size greater than
$2^{-m_n(1+\epsilon)}$ in $B(t_n,2^{-m_n\delta})$.

As we did before, up to extraction, we can e.g. assume that the first coordinate $\widetilde V^{B}$  of $V^B$  satisfies  $|\Delta \widetilde V_{t_n}^{B}| \ge 2^{-m_n}/3$ for all $n\geq 1$.

We then apply Lemma \ref{upbound} with $g(s)= \widetilde V_s^{B}$ and $g_n(s)=g(s)-\Delta \widetilde V_{t_n}^{B}\bbd{1}_{\{s\ge t_n\}}$,
with the choices $r_n=2^{-m_n\delta}$ and $\delta_n=m_n\delta 2^{-m_n(1+\e)}$.
It indeed holds true that $\lim_{n \to +\infty} \delta_n/|\Delta \widetilde V_{t_n}^{B}|=0$ and,
thanks to Lemma \ref{l2}-(1) (which is licit because $\delta/(1+\e)>\nu$) and the Borel-Cantelli Lemma,
we deduce that a.s., for all $n$ sufficiently large,
\[\osc_{B(t_n,r_n)}  (\widetilde V^{B}_s  )
\le \sup_{x,y\in[0,1], |x-y|\le 2^{-m_n\delta}}|V_x^{B,m_n(1+\e)} - V_y^{B,m_n(1+\e)}|\le m_n\delta 2^{-m_n(1+\e)}.\]
We conclude from Lemma \ref{upbound} that
\[h_{X^B}(t)\le\liminf_{n}
\frac{\log\Big(r_n|\Delta \widetilde V_{t_n}^{B}|\Big)}{\log(|t_n-t|+r_n)}\le \liminf_{n}
\frac{\log\Big(2^{-m_n(1+\delta)}/3\Big)}{\log(2.2^{-m_n/(h+\e)})}
=(1+\delta)(h+\e).\]
We used that $r_n|\Delta \widetilde V_{t_n}^{B}|\geq (2^{-m_n}/3)2^{-m_n\delta}$ while $|t_n-t|+r_n\leq 2^{-m_n/(h+\e)}
+2^{-m_n\delta}\leq 2.2^{-m_n/(h+\e)}$.
Letting $\e\to 0$ (whence $\delta \to 1/h$), we conclude that
$h_{X^B}(t) \leq 1+h$ as desired.
\epf

Before computing the dimension of $E^{osc}_{V^B}(h)$ when $h\in [1/(2\nu),1/\nu)$,
we need to count those jump times that are very close to each other.

\blem\label{numb}
For $\epsilon>0$ and $m>0$, denote by $0<T_1^{\e,m}<\dots<T_{K_{\e,m}}^{\e,m}<1$ the successive instants of jumps
of $V^B$ with size greater than $2^{-m(1+\epsilon)}$. For $\delta>0$, we introduce
\[
N_m^{\delta,\e}= \sum_{i=1}^{K_{\e,m}} \bbd{1}_{\{T_i^{\e,m}-T_{i-1}^{\e,m} \leq 2^{-m \delta}\}}
\]
with the convention that $T_0^{\e,m}=0$. For any fixed $\e>0$ and $\delta>0$, there a.s. exists $M>0$
such that for all $m>M$,
\[N_m^{\delta,\e}\le 2^{-m\delta+2m\nu(1+2\epsilon)}.\]
\elem
\bpf
Recalling Lemma \ref{le3}, we see that $\{T_1^{\e,m},\dots,T_{K_{\e,m}}^{\e,m}\}
\subset \{S_1^{\e,m},\dots,S_{L_{\e,m}}^{\e,m} \}$, where
$0<S^{\e,m}_1<\dots<S^{\e,m}_{L_{\e,m}}$ are the successive instants of jump of the counting process
$R^{m(1+\e)}_{[0,t]}$. Consequently,
$$N_m^{\delta,\e}\leq \tilde N_m^{\delta,\e}:= \sum_{i=1}^{L_{\e,m}} \bbd{1}_{\{S_i^{\e,m}-S_{i-1}^{\e,m} \leq 2^{-m \delta}\}}.$$
By Remark \ref{rem1}, we know that $R^{m(1+\e)}_{[0,t]}$ is an inhomogeneous Poisson process
with intensity bounded by $C_B 2^{m(1+\e)\nu}$. Consequently,
\[\bb{P}\Big[L_{\e,m}\ge 2^{m\nu(1+2\epsilon)}\Big]\le
2^{-m\nu(1+2\epsilon)}C_B 2^{m(1+\e)\nu} \leq C_B2^{-m\nu\epsilon}.\]
Hence, applying the Borel-Cantelli lemma, we know that almost surely, there exists $M'\ge 1$
such that for all $m\ge M'$,
\[L_{\e,m}\le 2^{m\nu(1+2\epsilon)}
\quad\hbox{and thus}\quad N_m^{\delta,\e} \leq \sum_{i=1}^{2^{m\nu(1+2\epsilon)}} \bbd{1}_{\{S_i^{\e,m}-S_{i-1}^{\e,m} \leq 2^{-m \delta}\}}.
\]
But for all $i\geq 1$, $S_i^{\e,m}-S_{i-1}^{\e,m}$ is bounded from above by an exponential
random variable with parameter $C_B 2^{m(1+\e)\nu}$, so that $\bb{P}(S_i^{\e,m}-S_{i-1}^{\e,m}\leq 2^{-m\delta})
\leq 1-\exp(-C_B 2^{m(1+\e)\nu} 2^{-m\delta})\leq C_B 2^{m(1+\e)\nu-m\delta}$ and thus
\begin{align*}
\bb{P}\Big(\sum_{i=1}^{2^{m\nu(1+2\epsilon)}} \bbd{1}_{\{S_i^{\e,m}-S_{i-1}^{\e,m} \leq 2^{-m \delta}\}} \ge2^{-m\delta+2m\nu(1+2\epsilon)}\Big)
\le&  \ 2^{m\delta-2m\nu(1+2\epsilon)}2^{m\nu(1+2\epsilon)}  C_B 2^{m(1+\e)\nu-m\delta} \\
=& \  C_B 2^{-m\nu\e}.
\end{align*}
By the Borel-Cantelli lemma again, there exists a.s.  a constant $M''>0$ such that
for all $m\geq M''$,
$$\sum_{i=1}^{2^{m\nu(1+2\epsilon)}} \bbd{1}_{\{S_i^{\e,m}-S_{i-1}^{\e,m} \leq 2^{-m \delta}\}} \le 2^{-m\delta+2m\nu(1+2\epsilon)}.$$
As a conclusion, a.s.   we have $N_m^{\delta,\e} \le 2^{-m\delta+2m\nu(1+2\epsilon)}$ for all $m\geq M'\lor M''$.
Choosing $M=M'\lor M''$ completes the proof.
\epf

Now we treat the case where $h \in [1/(2\nu),1/\nu)$.

\bprop\label{dim}
Almost surely, for $h\in [1/(2\nu),1/\nu)$, $\dim_H\Big(E^{osc}_{V^B}(h)\Big) \le 2h\nu-1$.
\eprop
\bpf We divide the proof into several steps.

\smallskip

{\it Step 1.}
For any fixed $\epsilon>0$, $\delta\in(\nu,2\nu]$ and $m\ge 1$, we consider the sets
$$F_m(\delta,\epsilon)=\bigcup_{\{i:T_i^{\e,m}-T_{i-1}^{\e,m} \leq 2^{-m \delta}\}}
\left([T_{i-1}^{\e,m}-2^{-m\delta},~T_{i-1}^{\e,m}+2^{-m\delta}]
\cup[T_i^{\e,m}-2^{-m\delta},~T_i^{\e,m}+2^{-m\delta}]\right),$$
 where the family $T_{i}^{\e,m}$
has been introduced in Lemma \ref{numb}, and the associated limsup set
\[G(\delta,\epsilon)=\limsup_{m\rightarrow+\infty}F_m(\delta,\epsilon)
.\]

For every $n\geq 1$, $\bigcup_{m\ge n}F_m(\delta,\epsilon)$ forms a covering of $G(\delta,\epsilon)$ by sets of diameter less than $2^{-n\delta+2}$, and Lemma \ref{numb} allows to bound by above the cardinality of such sets. Hence, choosing  $s > \frac{ 2\nu(1 + 2\epsilon)}{\delta}-1$, a.s. for every $n$ large enough one has
\[\ca{H}_{2^{-n\delta+2}}^s(G(\delta,\epsilon))\le \sum_{m\ge n}2^{-m \delta s + 2s} N_m^{\delta,\epsilon} \le \sum_{m\ge n} 2^{2s}2^{-m(s+1)\delta+2m\nu(1+2\epsilon)}.\]
We deduce that  $\lim_{n\rightarrow+\infty}\ca{H}_{2^{-n\delta+2}}^s(G(\delta,\epsilon))=0$, hence $\ca{H} ^s(G(\delta,\epsilon))=0$. Therefore, $\dim_H \Big(G(\delta,\epsilon)\Big)\le\frac{2\nu(1+2\epsilon)}{\delta}-1$.

\smallskip

{\it Step 2.} Here we fix $h\in [1/(2\nu),1/\nu)$, we consider
$\e>0$ such that $1/[(h+\e)(1+\e)]>\nu$, we set $\delta_\e=1/(h+\e)$
and we prove that $E^{osc}_{V^B}(h) \subset G(\delta_\e,\e)$.

We consider
$t\in E_{V^B}(h)\setminus G(\delta_\e,\e)$ and we show that $h_{X^B}(t)= 1+h$,  which will imply indeed that
$t\in E_{V^B}^{cusp}(h)$.

Since $t\notin G(\delta_\e,\e)$, there exists $N\ge 1$ such that for all $m \ge N$,
$t\notin F_m(\delta_\e,\e)$. Moreover, for any $0<\eta \le \e$,
since $t\in E_{V^B}(h)$, by Remark \ref{nostar}, we know that
$t \in A_{\delta_\eta}$ (because $\delta_\eta=1/(h+\eta)<1/h$), so that
for all $n \ge 1$, there exist $m_n \geq n$ and $t_n\in B(t, 2^{-m_n\delta_\eta})$
such that $|\Delta V_{t_n}^B|\ge 2^{-m_n}$. Observing that $F_m(\delta_\eta,\eta)\subset F_m(\delta_\e,\e)$ since $0<\eta \le \e$ and $\delta_\eta\ge \delta_\e$.
Hence $t\notin F_{m_n}(\delta_\eta,\eta)$ (for all $n$ large enough), whence, there is also no other jump in $B(t, 2^{-m_n\delta_\eta})$
with size greater than $2^{-m_n(1+\eta)}$.

As in the previous proofs, up to extraction, we deduce that $|\Delta \widetilde V_{t_n}^{B}| \ge 2^{-m_n}/3$
for all $n$, where $ \widetilde V_{ }^{B}$ is one of the three coordinates of $V^B$. Since $V^B$ (and so $ \widetilde V^{B}$
) has
no jump with size greater than $2^{-m_n(1+\eta)}$ in $B(t_n,2^{-m_n\delta_\eta})$,
we may use Lemma \ref{l2}-(1) (because $\delta_\eta/(1+\eta)=\frac{1}{(h+\eta)(1+\eta)}
\ge\frac{1}{(h+\e)(1+\e)}>\nu$)
and the Borel-Cantelli Lemma, we deduce that a.s. for all $n$ sufficiently large,
setting $r_n=2^{-m_n\delta_\eta}$,
\[\osc_{B(t_n,r_n)} (\widetilde V^{B})
\leq 2 \times
\sup_{x,y\in[0,1], |x-y|\le 2^{-m_n\delta_\eta}}|V_x^{B,m_n(1+\eta)} - V_y^{B,m_n(1+\eta)}|
\le 2m_n\delta_\eta 2^{-m_n(1+\eta)}.\]
Moreover,
\[\lim_{n\rightarrow+\infty}\frac{2m_n \delta_\eta 2^{-m_n(1+\eta)}}{|\Delta \widetilde V_{t_n}^{B}|}\le
\lim_{n\rightarrow+\infty}\frac{2m_n \delta_\eta 2^{-m_n(1+\eta)}}{2^{-m_n}/3}=0.\]
Applying Lemma \ref{upbound} with $g=\widetilde V^{B}$, $r_n= 2^{-m_n\delta_\eta}$ and
$\delta_n=2m_n \delta_\eta 2^{-m_n(1+\eta)}$, we obtain

\beq\label{final}
h_{X^B}(t) \le \liminf_{n\rightarrow+\infty}
\frac{\log\Big(r_n |\Delta \widetilde V_{t_n}^{B}|\Big)}{\log(r_n+|t_n-t|)}
\le\liminf_{n\rightarrow+\infty}
\frac{\log\Big(2^{-m_n(1+\delta_\eta)}/3\Big)}{\log (2.2^{-m_n\delta_\eta})}
=\frac{1+\delta_\eta}{\delta_\eta}=1+h+\eta
\eeq
because $r_n |\Delta \widetilde V_{t_n}^{B}|\geq 2^{-m_n(1+\delta_\eta)}/3$ and
$r_n+|t_n-t|\leq 2.2^{-m_n\delta_\eta}$. Since (\ref{final}) is satisfied for any $0<\eta \le \e$, then a.s. $h_{X^B}(t) \le 1+h$. That is, $E^{osc}_{V^B}(h) \subset G(\delta_\e,\e )$.

{\it Step 3.}
From step 2 we deduce that $E^{osc}_{V^B}(h)\subset \bigcap_{\e \downarrow 0}G(\delta_\e,\e)$. Hence,
\[\dim_H\Big(E^{osc}_{V^B}(h)\Big) \le \dim_H \Big(\bigcap_{\e \downarrow 0}G(\delta_\e,\e)\Big)
=\inf_{\e \downarrow 0} \Big(2\nu(1+2\epsilon)(h+\e)-1 \Big)=2h\nu - 1.\]
This ends the proof.
\epf

\subsection{Conclusion}

\bpf[Proof of Proposition \ref{locpo}]
First, we now from Proposition \ref{loc} that $E_{V^B}(h)=\emptyset$ for $h>1/\nu$, so that
obviously $E_{V^B}^{osc}(h)=\emptyset$. If now $h=1/\nu$, then we deduce from
Proposition \ref{uniform} that $E_{V^B}^{osc}(h)=\emptyset$, simply because
a.s., for all $t\in [0,1]$, $h_{X^B}(t) \leq 1+1/\nu$.

As shown in Lemma \ref{rela}, we also know that $E_{V^B}^{osc}(h)=\emptyset$ for all $h \in [0,1/(2\nu))$
and as seen in Proposition \ref{dim},  $\dim_H (E_{V^B}^{osc}(h))\leq 2 h \nu - 1$ for all $h \in [1/(2\nu),1/\nu)$.

It remains to verify that for all $h\in [0,1/\nu]$, $\dim_H (E_{V^B}^{cusp}(h))=h \nu$ .
If $h \in [0,1/(2\nu))$ or $h=1/\nu$, it is obvious because $E_{V^B}^{osc}(h)=\emptyset$
and by Proposition \ref{loc}. If next $h \in [1/(2\nu),1/\nu)$, it follows from
the fact that $E_{V^B}^{cusp}(h)=E_{V^B}(h) \setminus E_{V^B}^{osc}(h)$ with
$\dim_H (E_{V^B}(h)) = h \nu$ (by Proposition \ref{loc}) and $\dim_H (E_{V^B}^{osc}(h))\leq 2 h \nu - 1<h \nu$.
\epf

Finally, we verify that Theorems \ref{m} and \ref{pos2} imply Theorem \ref{pos}.

\begin{proof}[Proof of Theorem \ref{pos}]
For any $h\in [1,1+1/\nu]$, we have $E_X(h) \supset E_V^{cusp}(h-1)$, whence
$\dim_H(E_X(h)) \geq  \dim_H(E_V^{cusp}(h-1)) = (h-1)\nu$ by Theorem \ref{pos2}.

Next we obviously have a.s., for all $t\in [0,1]$,
\begin{equation}\label{axa}
h_{X}(t) \geq h_V(t)+1,
\end{equation}
whence $E_X(h) \subset \bigcup_{h'\leq h-1} E_V(h')$. We thus infer from Theorem \ref{m} that
$E_X(h)=\emptyset$ when  $h<1$. But when  $h\in [1,1+1/\nu]$, recalling Proposition \ref{enplus} and the fact that $A_\delta^*$ is decreasing with $\delta$, we deduce that $\bigcup_{h'\leq h-1} E_V(h') \subset \bigcup_{h'\leq h-1}\bigcap_{\delta\in (0,1/h')} A_\delta^* \subset \bigcap_{\delta< h-1}A_\delta^*$. Whence we derive $\dim_H(E_X(h)) \leq (h-1)\nu$ from Proposition \ref{hd}.

It only remains to verify that $E_X(h)=\emptyset$ when  $h>1+1/\nu$. But in such a case, we know
from Proposition \ref{uniform} that $E_{X^B}(h)=\emptyset$, whence $E_X(h)=\bigcup_{B\ge 1}^{+\infty}E_{X^B}(h)=\emptyset$.
\end{proof}

\section*{Acknowledgement}
\quad~~
I would like to thank greatly my supervisors  Prof. Nicolas Fournier and Stephane Seuret for bringing me into this interesting subject with continuous support and encouragement on my study and research. I also wish to thank Prof. Stephane Jaffard for helpful discussions and kind help.

\bibliographystyle{amsplain}
\bibliography{ref}

\providecommand{\bysame}{\leavevmode\hbox to3em{\hrulefill}\thinspace}
\providecommand{\MR}{\relax\ifhmode\unskip\space\fi MR }
\providecommand{\MRhref}[2]{%
  \href{http://www.ams.org/mathscinet-getitem?mr=#1}{#2}
}
\providecommand{\href}[2]{#2}
\begin{thebibliography}{10}

\bibitem{alexandre2009review}
R.~Alexandre, \emph{A review of {Boltzmann} equation with singular kernels},
  Kinet. Relat. Models, \textbf{2} (2009), no.~4, 551--646.

\bibitem{arneodo1998singularity}
A.~Arneodo, E.~Bacry, S.~Jaffard, J.F. Muzy, \emph{Singularity spectrum of
  multifractal functions involving oscillating singularities}, J. Fourier Anal.
  Appl. \textbf{4} (1998), no.~2, 159--174.

\bibitem{balancca2014fine}
P.~Balan{\c{c}}a, \emph{Fine regularity of {L{\'e}vy} processes and linear
  (multi) fractional stable motion}, Electron. J. Probab \textbf{19} (2014),
  no.~101, 1--37.

\bibitem{barral2010pure}
J.~Barral, N.~Fournier, S.~Jaffard, S.~Seuret, \emph{A pure jump {Markov}
  process with a random singularity spectrum}, Ann. Probab. \textbf{38} (2010),
  no.~5, 1924--1946.

\bibitem{beresnevich2006mass}
V.~Beresnevich, S.~Velani, \emph{A mass transference principle and the
  {Duffin-Schaeffer} conjecture for {Hausdorff} measures}, Ann. of Math. (2)
  \textbf{164} (2006), 971--992.

\bibitem{cercignani1988boltzmann}
C.~Cercignani, \emph{The {Boltzmann} equation and its applications}, vol. 67 of
  Applied Mathematical Sciences, Springer-Verlag, New York, 1988.

\bibitem{falconer2007fractal}
K.~Falconer, \emph{Fractal geometry: mathematical foundations and
  applications}, Wiley, 2007.

\bibitem{fournier2012finiteness}
N.~Fournier, \emph{Finiteness of entropy for the homogeneous {Boltzmann}
  equation with measure initial condition}, Ann. Appl. Probab., to appear
  \textbf{25} (2015), no.~2, 860--897.

\bibitem{fournier2001markov}
N.~Fournier, S.~M{\'e}l{\'e}ard, \emph{A {Markov} process associated with a
  {Boltzmann} equation without cutoff and for non-{Maxwell} molecules}, J.
  Stat. Phys. \textbf{104} (2001), no.~1-2, 359--385.

\bibitem{fournier2002stochastic}
N.~Fournier, S.~M{\'e}l{\'e}ard, \emph{A stochastic particle numerical method for {3D} {Boltzmann}
  equations without cutoff}, Math. Comp. \textbf{71} (2002), no.~238, 583--604.

\bibitem{fm09}
N.~Fournier, C.~Mouhot, \emph{On the well-posedness of the spatially
  homogeneous {Boltzmann} equation with a moderate angular singularity}, Comm.
  Math. Phys. \textbf{289} (2009), no.~3, 803--824.

\bibitem{jacod1987limit}
J.~Jacod, A.N. Shiryaev, \emph{Limit theorems for stochastic processes},
  vol. 288, Springer-Verlag, 1987.

\bibitem{jaffard1999multifractal}
S.~Jaffard, \emph{The multifractal nature of {L{\'e}vy} processes}, Probab.
  Theory Related Fields \textbf{114} (1999), no.~2, 207--227.

\bibitem{lu2012measure}
X.~Lu, C.~Mouhot, \emph{On measure solutions of the {Boltzmann} equation,
  part {I}: moment production and stability estimates}, J. Differential
  Equations \textbf{252} (2012), no.~4, 3305--3363.

\bibitem{orey1974often}
S.~Orey, S.J. Taylor, \emph{How often on a {Brownian} path does the law of
  iterated logarithm fail?}, Proc. London Math. Soc. (3) \textbf{3} (1974),
  no.~1, 174--192.

\bibitem{perkins1983hausdorff}
E.~Perkins, \emph{On the {Hausdorff} dimension of the {Brownian} slow points},
  Z. Wahrsch. Verw. Gebiete \textbf{64} (1983), no.~3, 369--399.

\bibitem{shepp1972covering}
L.A. Shepp, \emph{Covering the line with random intervals}, Probab. Theory
  Related Fields \textbf{23} (1972), no.~3, 163--170.

\bibitem{tanaka1978probabilistic}
H.~Tanaka, \emph{Probabilistic treatment of the {Boltzmann} equation of
  {Maxwellian} molecules}, Z. Wahrsch. Verw. Gebiete \textbf{46} (1978), no.~1,
  67--105.

\bibitem{villani1998new}
C.~Villani, \emph{On a new class of weak solutions to the spatially homogeneous
  {Boltzmann} and {Landau} equations}, Arch. Ration. Mech. Anal., \textbf{143}
  (1998), no.~3, 273--307.

\bibitem{villani2002review}
C.~Villani, \emph{A review of mathematical topics in collisional kinetic theory},
  In Handbook of mathematical fluid dynamics \textbf{1} (2002), 71--305.

\end{thebibliography}

\vskip 3\baselineskip

Liping Xu\\
Laboratoire de Probabilit\'es et Mod\`eles Al\'eatoires, UMR 7599, Universit\'e Pierre-et-Marie Curie,
Case 188, 4 place Jussieu, F-75252 Paris Cedex 5, France.\\
E-mail: \texttt{xuliping.p6@gmail.com}.

\end{document}